\documentclass{ws-ijbc}
\usepackage{ws-rotating}     
\usepackage[english]{babel}
\usepackage{wasysym}
\usepackage[usenames]{color}
\usepackage{graphicx}

\begin{document}
\renewcommand{\figurename}{Fig.}
\markboth{Albert D. Morozov and Olga S. Kostromina}{On periodic
perturbations of asymmetric Duffing--Van-der-Pol equation}
\title{\bfseries{ON PERIODIC PERTURBATIONS OF ASYMMETRIC DUFFING--VAN-DER-POL EQUATION}}
\author{Albert D. Morozov  and Olga S. Kostromina}
\address{\footnotesize{\it{Department of Mechanics and Mathematics,
     Lobachevsky State University of Nizhny Novgorod,
\\ Nizhny Novgorod,  603950, Russia \\
  Correspondence should be addressed to Albert D. Morozov, morozov@mm.unn.ru}}}

\maketitle
\hyphenpenalty=10000 
\begin{abstract}
Time-periodic perturbations of an asymmetric Duffing--Van-der-Pol equation close to an integrable equation with a homoclinic ``figure-eight'' of a saddle are considered. The behavior of solutions outside the neighborhood of ``figure-eight'' is studied analytically. The problem of limit cycles for an autonomous equation is solved and resonance zones for a nonautonomous equation are analyzed. The behavior of the separatrices of a fixed saddle point of the Poincar\'{e} map in the small neighborhood of the unperturbed ``figure-eight'' is ascertained. The results obtained are illustrated by numerical computations.
\end{abstract}


\section{Introduction}
\begin{figure}[htb]
\begin{center}
\includegraphics[scale=0.36]{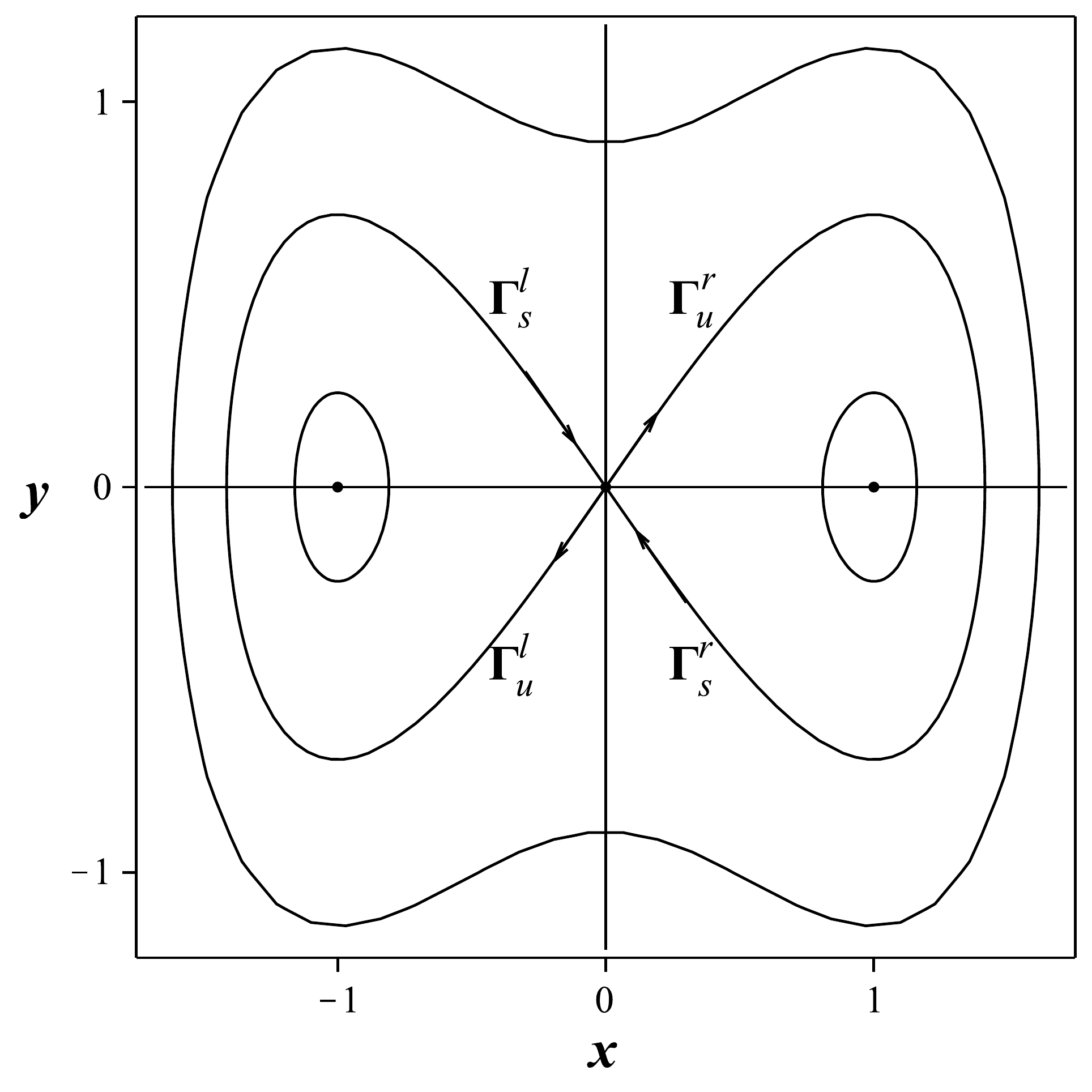}
\caption{Phase portrait of unperturbed equation.}\label{fig1}
\end{center}
\end{figure}

The theory of time-periodic systems close to two-dimensional nonlinear Hamiltonian systems has been greatly advanced by now (see, e.g., \cite{GH}, \cite{MorSh1983}, \cite{Wig}, \cite{Mor1998}). However, many problems remain unsolved, and new examples should be addressed. In this paper, we will consider one such example -- an asymmetric variant of the classical
Duffing--Van-der-Pol equation:
\begin{equation}\label{eq0}
\ddot{x}+\alpha x+\beta x^3=\varepsilon [(\gamma_1+\gamma_2
x+\gamma_3 x^2)\dot{x}+\gamma_4\sin{\gamma_5t}],
\end{equation}
where $\alpha , \beta, \gamma_1\div \gamma_5$ are parameters, and $\varepsilon$ is a small positive parameter. It is always possible
to set $\alpha =\pm 1, \beta =\pm 1$ in Eq.~ (\ref{eq0}). Naturally, the case $\alpha=\beta =-1$ is of no interest to us. In the cases $\alpha =1, \beta =\pm 1$,  the asymmetric perturbation term $\gamma_2 x\dot{x}$  does not play a significant role, it is essential only in the case $\alpha =-1, \beta=1$. Equation (\ref{eq0}) for $\alpha =1, \beta =\pm 1$ was studied in ample
detail (see, e.g., \cite{Mor1973}, \cite{MorSh1975}, \cite{Mor1976}, \cite{Mor1993}, \cite {Mor1998}). Therefore, we will address the case $\alpha =-1, \beta=1$.  The phase plane of an unperturbed equation has two saddle separatrix loops $O(0,0)$ forming ``figure-eight'' (Fig.~\ref{fig1}). Of the three parameters $\gamma_1, \gamma_2, \gamma_3$  one may be excluded to yield the following equation
\begin{equation}\label{eq1}
\ddot{x}-x+x^3=\varepsilon [(p_1+p_2x-x^2)\dot{x}+p_3\sin{p_4t}],
\end{equation}
where  $p_1\div p_4$ are parameters\footnote{The equation with
parametric perturbation was considered in \cite{LHRK}.}.
   An analysis of Eq.~(\ref{eq1}) implies the solution of the following tasks:
    1) for the autonomous equation ($p_3=0$) -- partition the plane
    of the parameters $(p_1,p_2)$ into regions with different topological structures and specify the structures;
    2) for the nonautonomous equation ($p_3\neq 0$) -- determine possible structures of resonance zones outside the neighborhood of ``figure-eight'' and the conditions of existence of structurally stable and unstable homoclinic Poincar\'{e} structures in the neighborhood of ``figure-eight''. The existence of a homoclinic structure specifies complicated behavior of solutions, in other words, it leads to chaos. Bifurcations in the neighborhood of ``figure-eight'' at a nonzero saddle value of the unperturbed autonomous system were recently considered in \cite{GSV}.

    The Duffing--Van-der-Pol equation is widely used in the theory of oscillations (see, e.g., \cite{Mor1973}, \cite{GH}, \cite{Mor1998}).
    Along with numerous applied problems in which there arises Eq.~(\ref{eq1}), we can mention a purely mathematical problem of vector field bifurcations on a plane that are invariant to the turn of angle $\pi $ \cite{Arnold}.
    In this problem, in Eq.~(\ref{eq0}) we have $\gamma_2=\gamma_4=0$ and the coefficients of the linear terms $\alpha x+\varepsilon \gamma_1 \dot{x}$, unlike our case, are the parameters of deformation. Note also the work \cite{Bautin} (as well as \cite{BautinLeontovich}) where an autonomous system with cubic nonlinearity without small parameter describing an electric circuit with tunnel diode was considered.  Possible local
bifurcations were determined and phase portraits were constructed
to an accuracy of an even number of limit cycles.

The presence of the term $p_2x\dot{x}$ in Eq. (\ref{eq1}) greatly
complicates the problem: there may exist in the autonomous equation ($p_3=0$) two limit cycles enclosing any of the equilibrium states $O_{\pm}(\pm 1,0)$, or a ``big'' separatrix loop enclosing the equilibrium states $O_{\pm}(\pm 1,0)$ that is absent in the unperturbed equation \cite{MorFed}, \cite{KM}.

Solution of problems 2) and 3) rests upon the solution of problem 1).
Therefore, we will start with the first problem.

\section{Investigation of an autonomous equation}

\subsection{Poincar\'{e}--Pontryagin generating functions}

The first integral of the unperturbed equation is $H(x,y)\equiv y^2/2-x^2/2+x^4/4=h$. The values $h\in (-0.25,0)$ correspond to
the domains  $G^{\pm}_1$ inside ``figure-eight'', and $h>0$ to the
domain $G_2$ outside ``figure-eight'' (Fig.~ \ref{fig1}). The value
$h=0$ corresponds to two symmetric saddle loops (``figure-eight'').

The main problem in studying Eq.~(\ref{eq1}) for $p_3=0$ is limit
cycles. Its solution results in finding real zeros of the
Poincar\'{e}--Pontryagin generating functions
$B(h)$~\cite{Mor1998}.

According to \cite{KM}, we have
\begin{equation}\label{eq6-1}
\begin{split}
B_1=B_1^{\pm}(\rho
(h))&=\frac{4}{30\pi(2-\rho)^{5/2}}\{2(5p_1-1)(\rho-1)(2-\rho)
{\bf K}(k)+\\
&+\left[5p_1(2-\rho)^2-4(\rho^2-\rho+1)\right]{\bf E}(k)\pm
\frac{15p_2\sqrt{2}}{16}\pi \rho^2\sqrt{2-\rho}\} \equiv \\
 &\equiv \frac{4}{30\pi(2-\rho)^{5/2}}B_{10}^{\pm}(\rho)
\end{split}
\end{equation}
for the domains $G_1^{\pm}$  and
\begin{equation}\label{eq7-1}
\begin{split}
B_2&=B_2(\rho
(h))=\frac{8}{30\pi(2\rho-1)^{5/2}}\left\{\left[5p_1(2\rho-1)(1-\rho)-
2(\rho-1)(2-\rho)\right]\right.{\bf K}(k)+\\
&+\left.\left[5p_1(2\rho-1)^2-4(\rho^2-\rho+1)\right]{\bf
E}(k)\right\}\equiv \frac{8}{30\pi(2\rho-1)^{5/2}}B_{20}(\rho)
 \end{split}
\end{equation}
for the domain $G_2$. Here, ${\bf K}(k), {\bf E}(k)$ are complete
elliptic integrals and $\rho=k^2$.  Note that $\rho$ is a more
convenient variable than $h$. In the formula  (\ref{eq6-1})  we
have $\rho(h)=2\sqrt{1+4h}/(1+\sqrt{1+4h})$ ($\rho\in (0,1)$), and
in (\ref{eq7-1}) we have $\rho(h)=(1+\sqrt{1+4h})/2\sqrt{1+4h}$
($\rho\in (1/2,1)$); $\rho =1 $ corresponds to ``figure-eight''. Note
that the function $B_2(\rho)$ does not depend on the parameter
$p_2$.

\subsection{Limit cycles}
Investigation of the functions $B_1^{\pm}(\rho ), B_2(\rho )$
gives the following results \cite{KM}.
\begin{theorem}\label{t1}
For sufficiently small  $\varepsilon $, the number of limit cycles
in each domain $G_1^{\pm}$  and $G_2$ of Eq.~(\ref{eq1}) does not
exceed two.
\end{theorem}
\begin{theorem}\label{t2} For sufficiently small $\varepsilon $,
the number of limit cycles of Eq.~(\ref{eq1}) does not exceed three.
\end{theorem}

The authors of \cite{KM} partitioned the parameter plane into 22
domains $D_m$, $m=1,\ldots,22$, gave the basic phase  portraits at
different values of parameters from those domains and found all
bifurcations.

\begin{figure}[htb]
\begin{center}
\begin{tabular}{cc}
\includegraphics[scale=0.7]{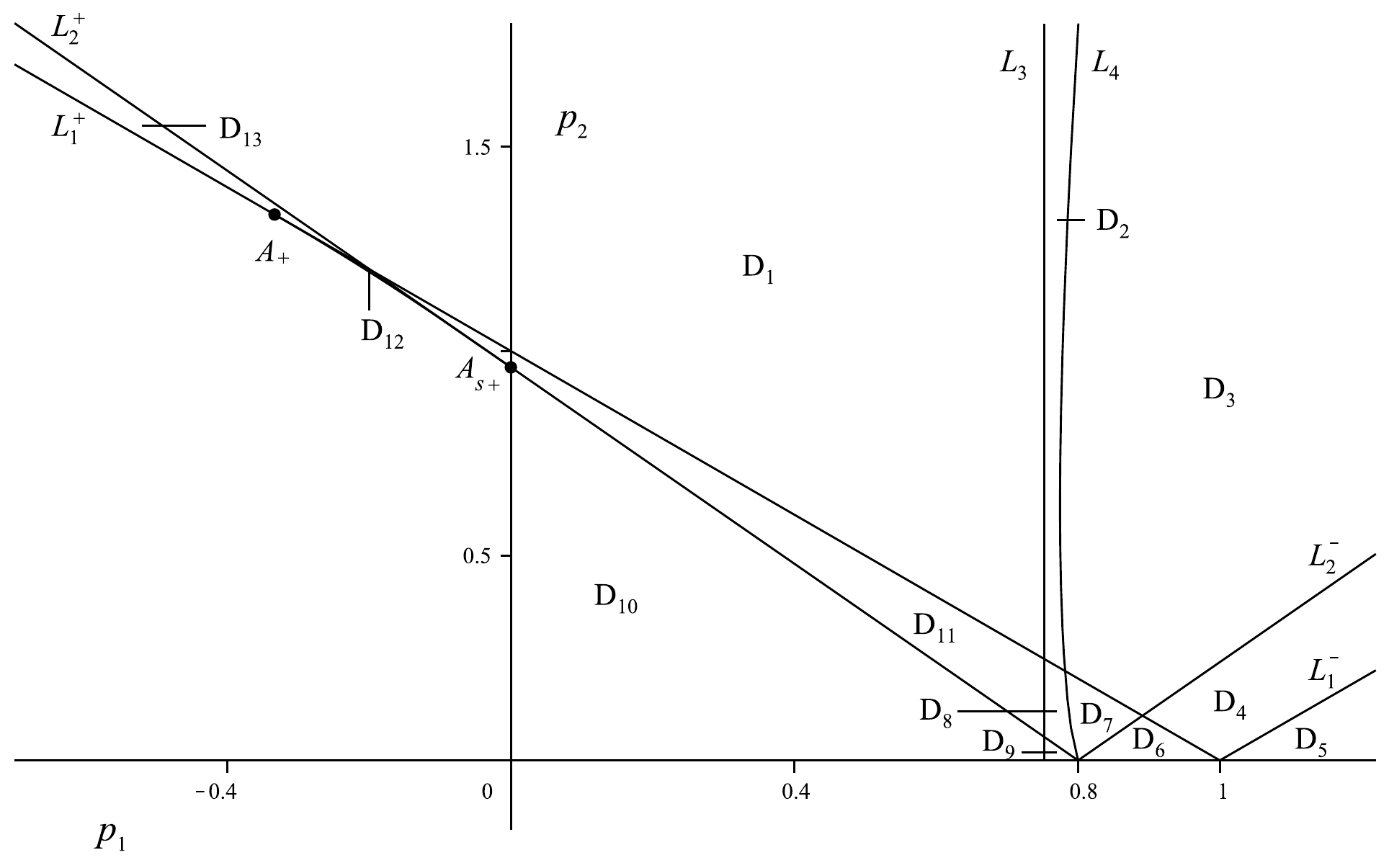}&
\end{tabular}
\caption{Partition of the plane of parameters $(p_1,p_2)$ into
domains with different phase portrait topology.} \label{fig2}
\end{center}
\end{figure}

\begin{figure}[htb]
  \begin{center}
    \begin{tabular}{c}
    \includegraphics[scale=0.4]{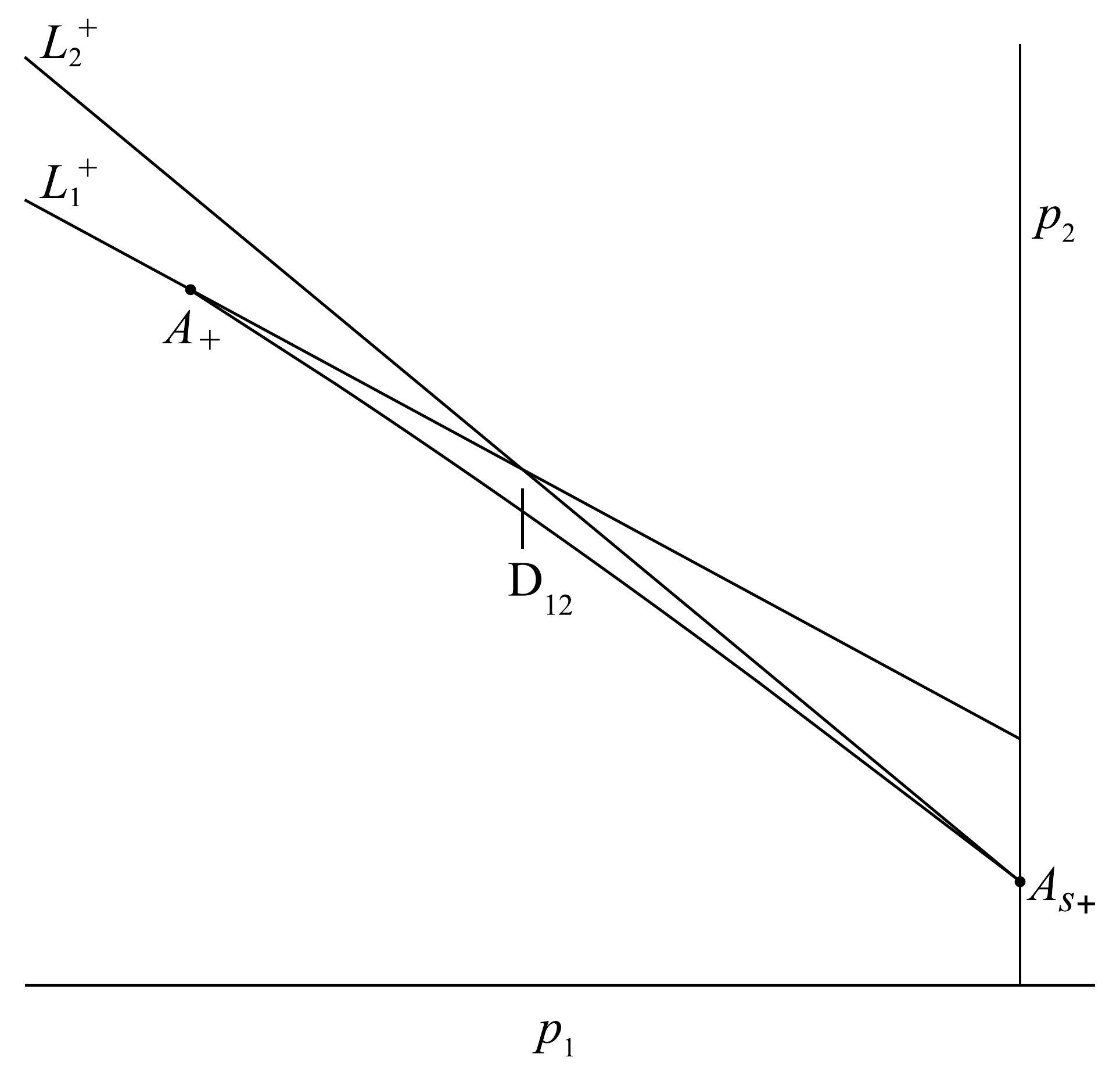}
    \end{tabular}
\end{center}
\caption{Magnified fragment of Fig.~\ref{fig2}.} \label{fig3}
\end{figure}

By virtue of the invariance of Eq.~(\ref{eq1}) to the change
$(p_2,x,y)\to (-p_2,-x,-y)$  the partitioning of the plane of the
parameters $(p_1,p_2)$ is symmetric to the $p_1$ axis. Therefore,
only the upper half plane $p_2\geq 0$ containing 13 domains $D_m$,
$m=1,\ldots,13$ is shown in Fig.~\ref{fig2}. A magnified fragment
of Fig.~\ref{fig2} is presented schematically in Fig.~\ref{fig3}.
By virtue of the symmetry, the phase portraits for the domain
$p_2<0$ are obtained by the turn of angle $\pi$ of the
corresponding phase portrait from the domain $p_2>0$.

Let us introduce the notation $(i,j,k)$ that means the existence
of $i$ limit cycles inside the right loop, $j$ inside the left
loop, and $k$ outside ``figure-eight''. It was proved that $D_1$ is of
the type $(0,0,0); D_2 - (0,0,2), D_3 - (0,0,1); D_4 - (0,1,1);
D_5 - (0,0,1); D_6 - (1,1,1); D_7 - (1,0,1); D_8 - (1,0,2); D_9
-(0,0,2); D_{10} - (0,0,0); D_{11} - (1,0,0); D_{12} - (2,0,0);
D_{13} - (1,0,0)$.

The notation of the bifurcation lines in Fig.~\ref{fig2}:

$L_1^{\pm}: p_1\pm  p_2-1=0$ -- straight lines at which the
autonomous equation has structurally unstable foci $O_{\pm}(\pm
1,0)$ in domains $G_1^{\pm}$, respectively.

 $A_{+} (-\frac {1} {3}, \frac {4} {3}) $ -- the point on the straight
 line $L_1^{+}$ from which the double cycle line originates. The double cycle line is plotted using the system $B_{10}^{+}(\rho , p_1,p_2)=0,
 [dB_{10}^{+}(\rho , p_1,p_2)/d\rho ]=0$, $\rho \in (0,1)$; the point  $A_{+}$ corresponds to $\rho =0$. This point is readily found by the power series expansion of the function $B_{10}^{+} (\rho )$ in the neighborhood of $\rho = 0$ in the case of a structurally unstable focus and zero first Lyapunov exponent (with the second Lyapunov exponent being nonzero).
 The extreme point $A_{s+}(0, 0.96)$ of the double cycle line corresponds to $\rho =1$ (see Fig.~\ref{fig3}). When the saddle value $\sigma _c=\varepsilon p_1$ vanishes to zero at $p_1=0$ the double cycle merges with the separatrix\footnote{In the lower half plane we have the points
 $A_{-} (-\frac {1} {3}, -\frac {4} {3})$ and $A_{s-} (0,-0.96) $, respectively.}.

$L_{2}^{\pm}: 5p_1\pm \frac{15\sqrt{2}\pi }{16}p_2-4=0$ -- straight lines on which the autonomous equation has separatrix loops (right, left) of saddle $O(0,0)$ in domains $G_1^{\pm}$. These straight lines are given by the Melnikov formula for autonomous systems.

$L_3$ -- double cycle line in domain $G_2$ corresponding to the
bifurcation value $p_1\approx 0.7523$.

$L_4$ -- line of the ``big'' loop of the separatrix of saddle
$O(0,0)$. This line was plotted numerically using the WInSet
software \cite{MD}.

\begin{figure}[htb]
\begin{center}
\begin{tabular}{ccc}
\includegraphics[scale=0.35]{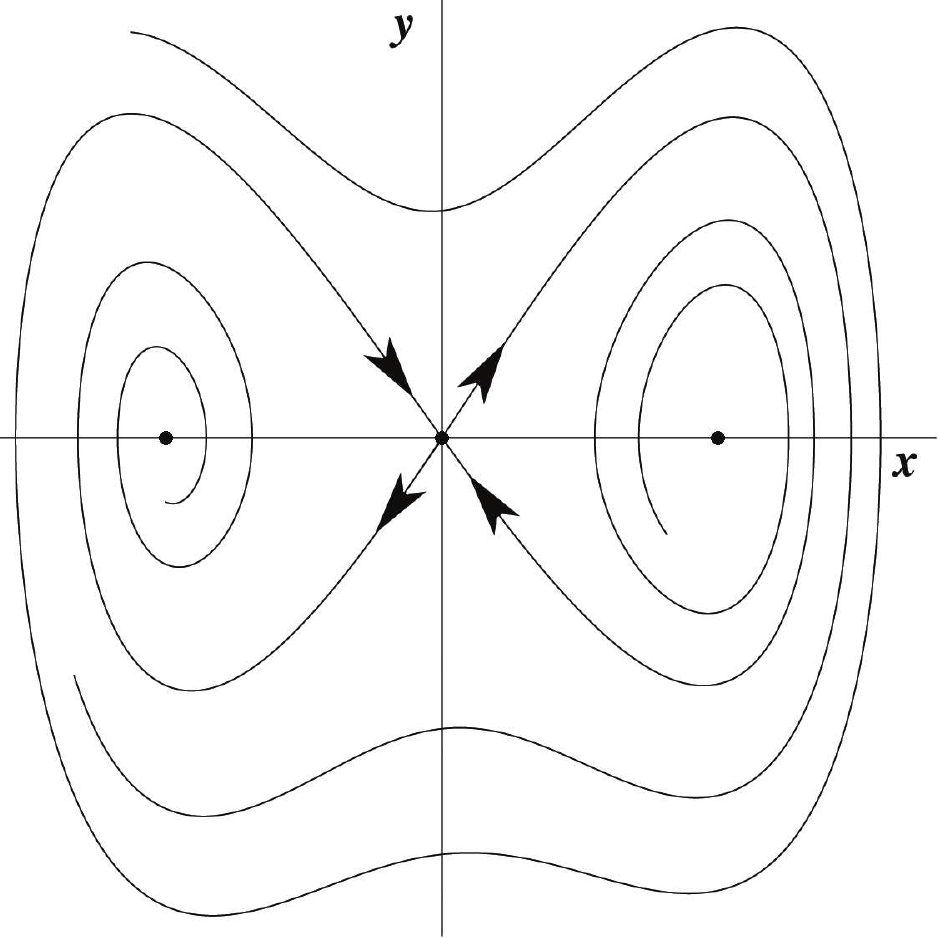}&
\includegraphics[scale=0.35]{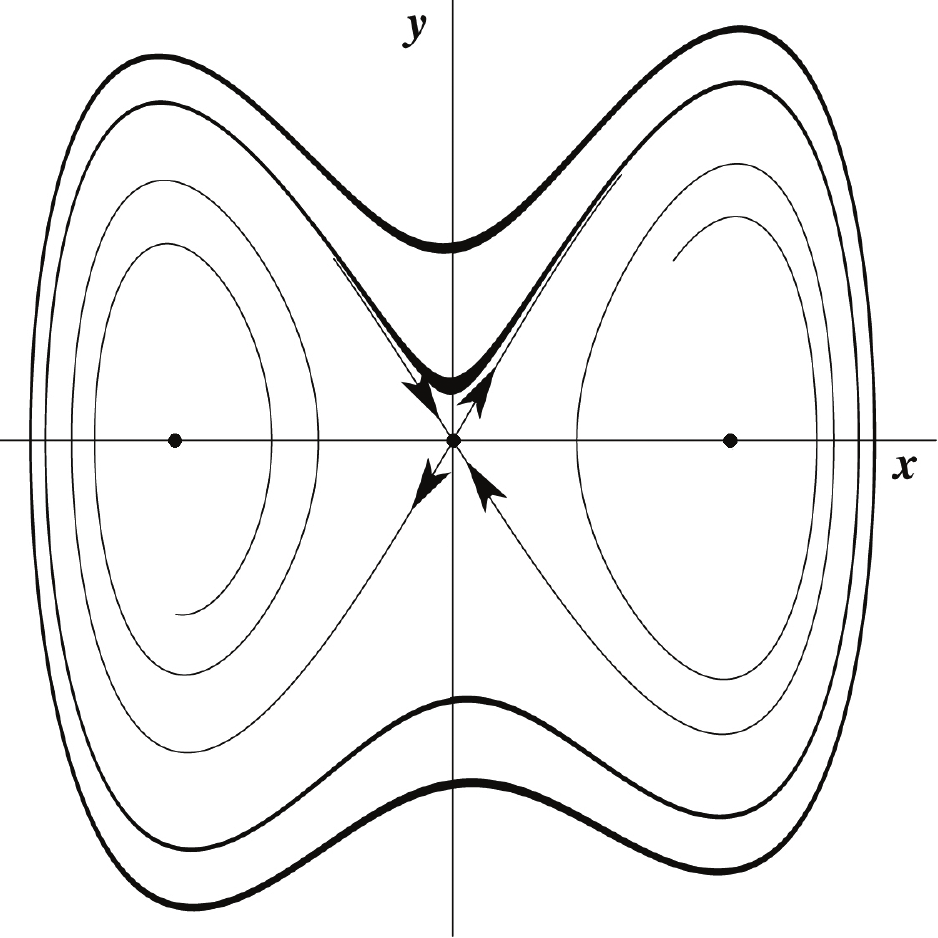}&
\includegraphics[scale=0.35]{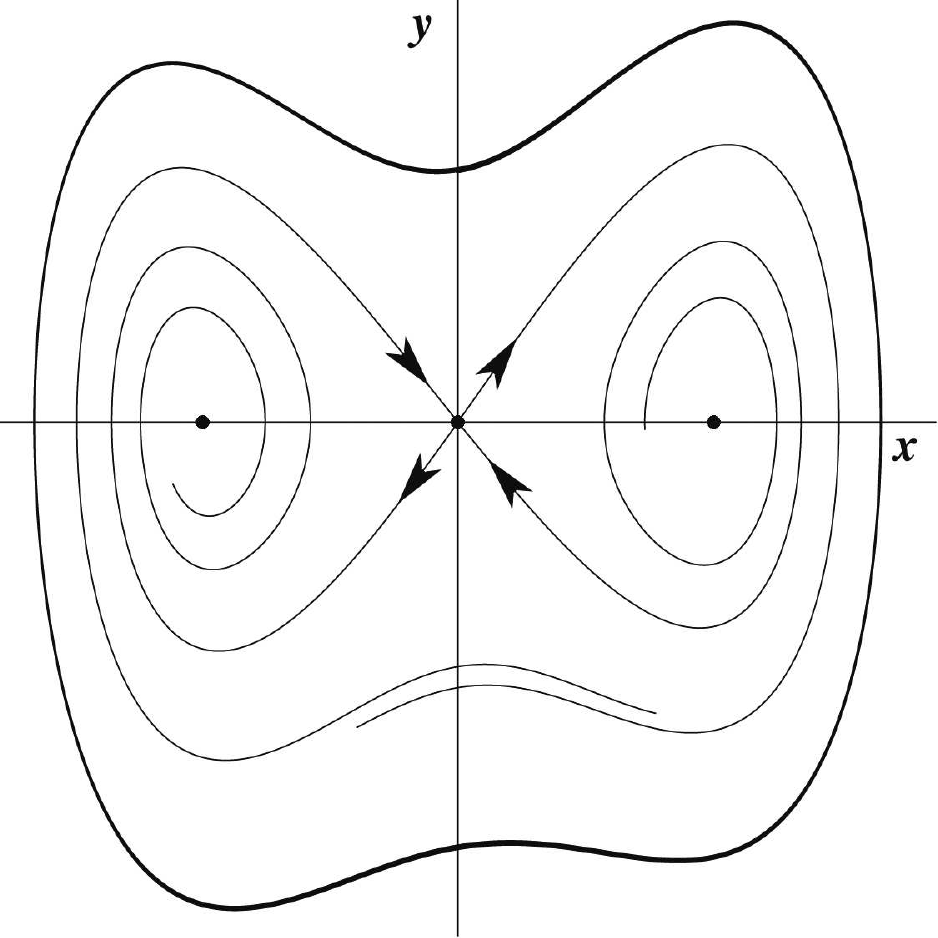}\\
\footnotesize{(a)}& \footnotesize{(b)} & \footnotesize{(c)}\\
\includegraphics[scale=0.35]{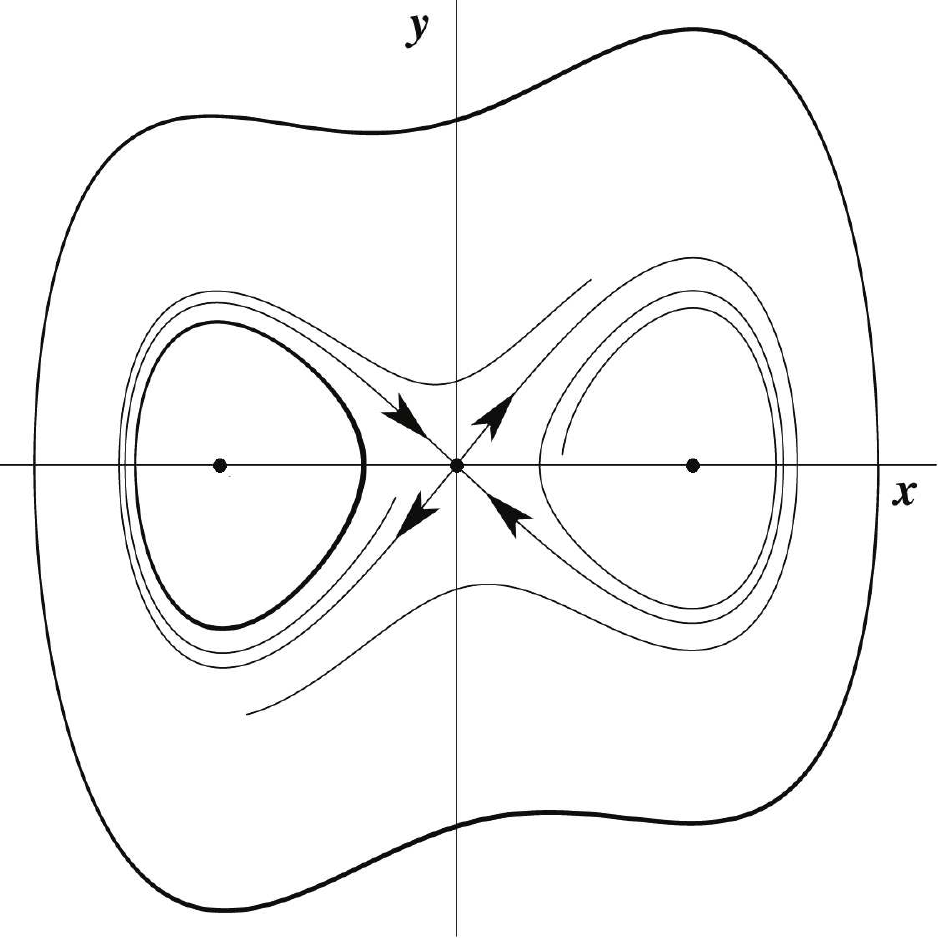}&
\includegraphics[scale=0.35]{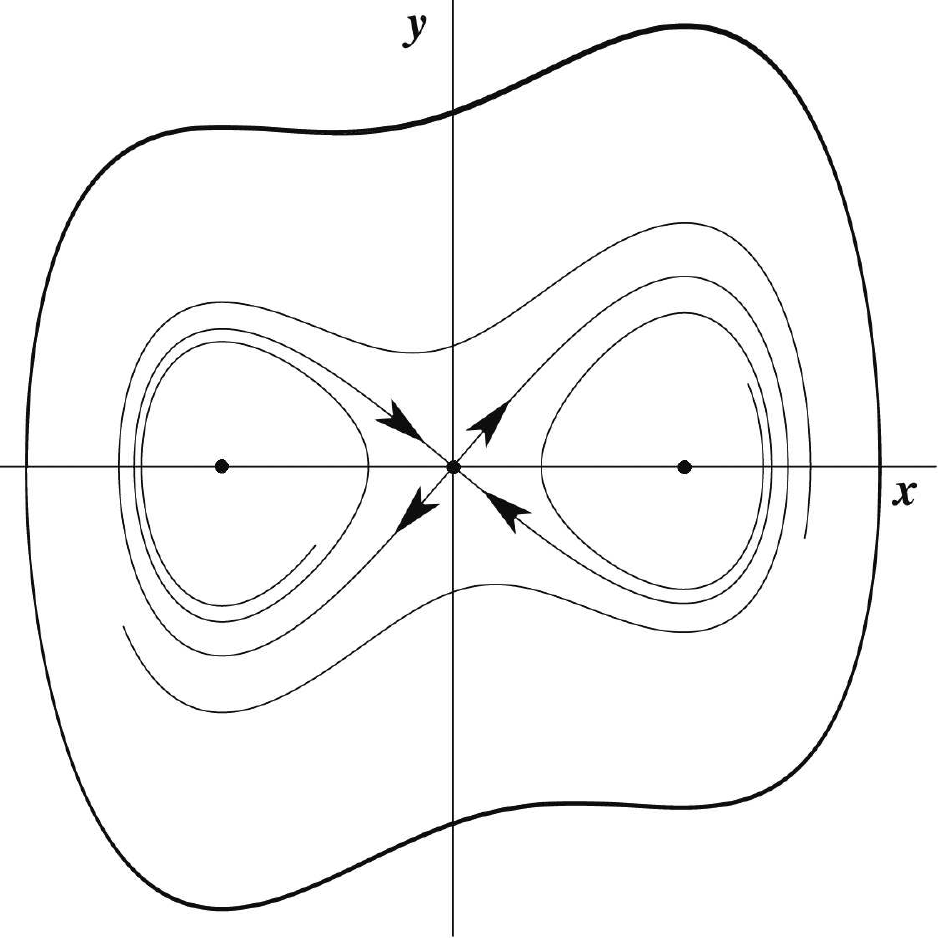}&
\includegraphics[scale=0.35]{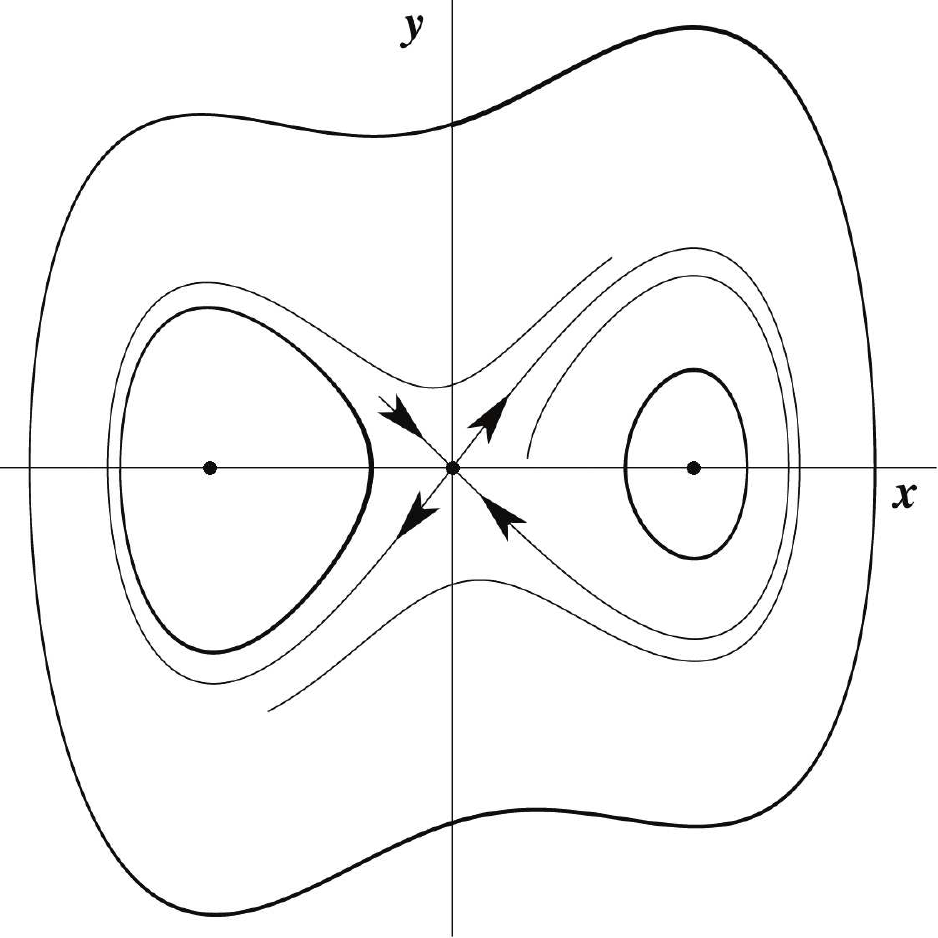}\\
\footnotesize{(d)}& \footnotesize{(e)} & \footnotesize{(f)}\\
\includegraphics[scale=0.35]{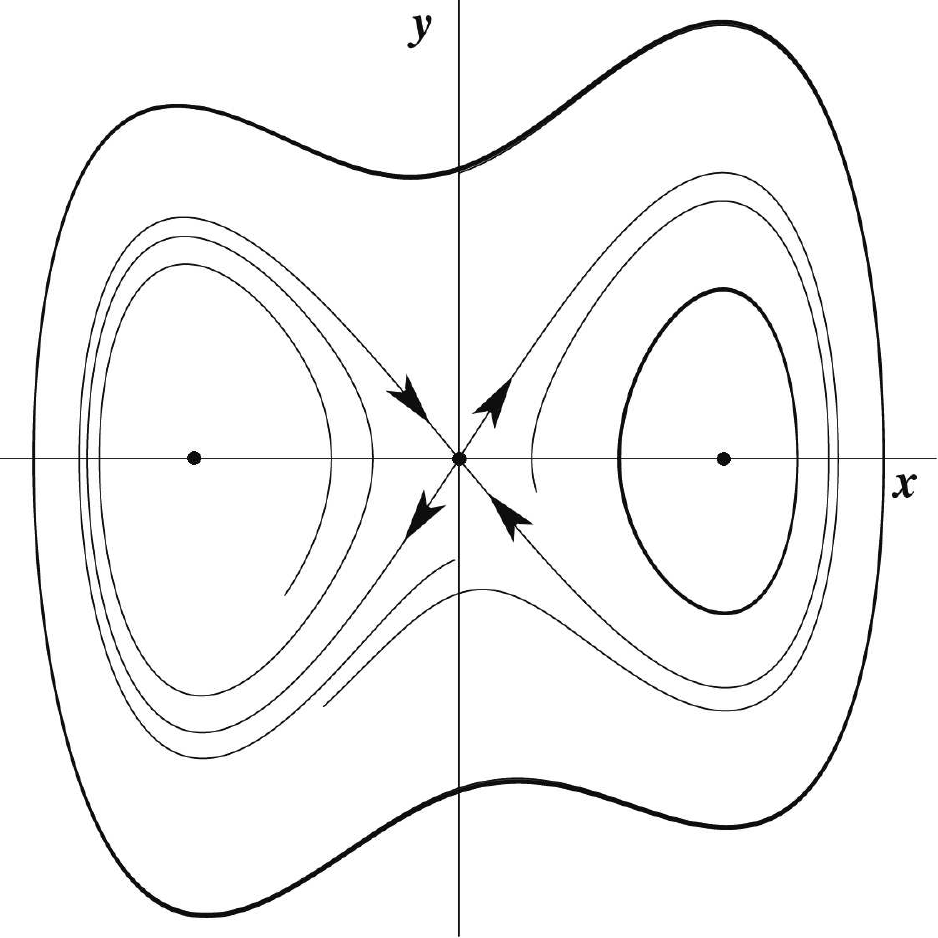}&
\includegraphics[scale=0.35]{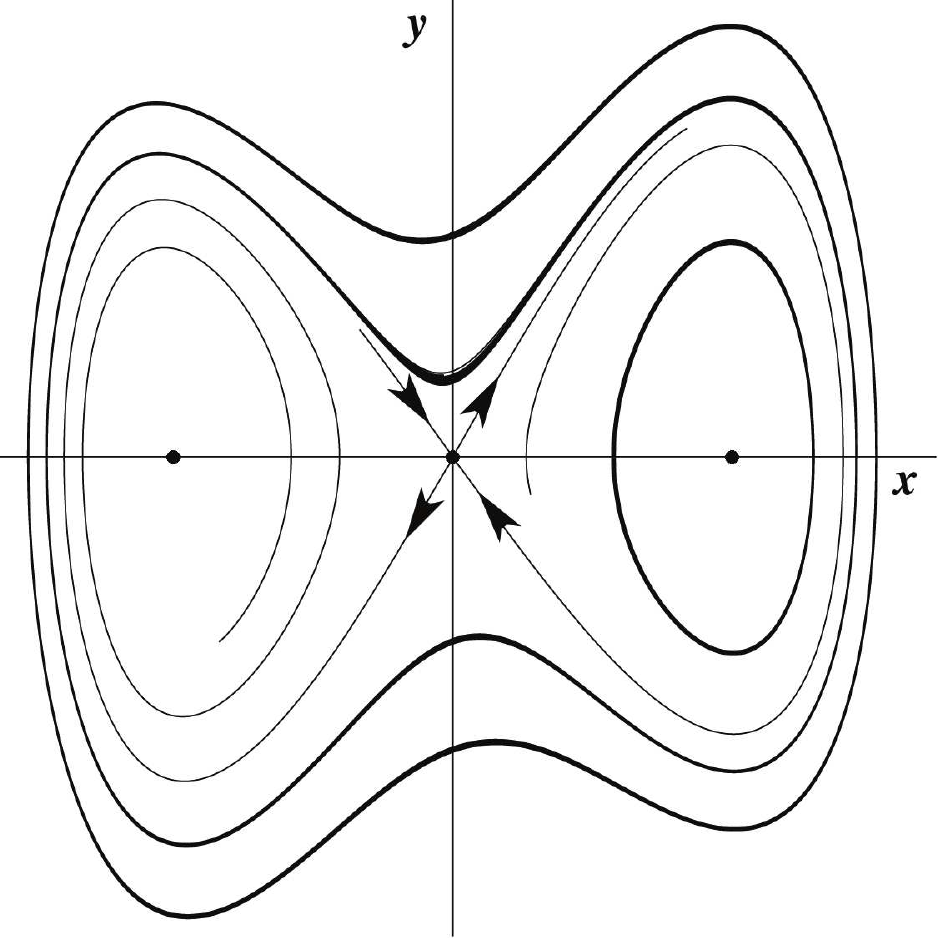}&
\includegraphics[scale=0.35]{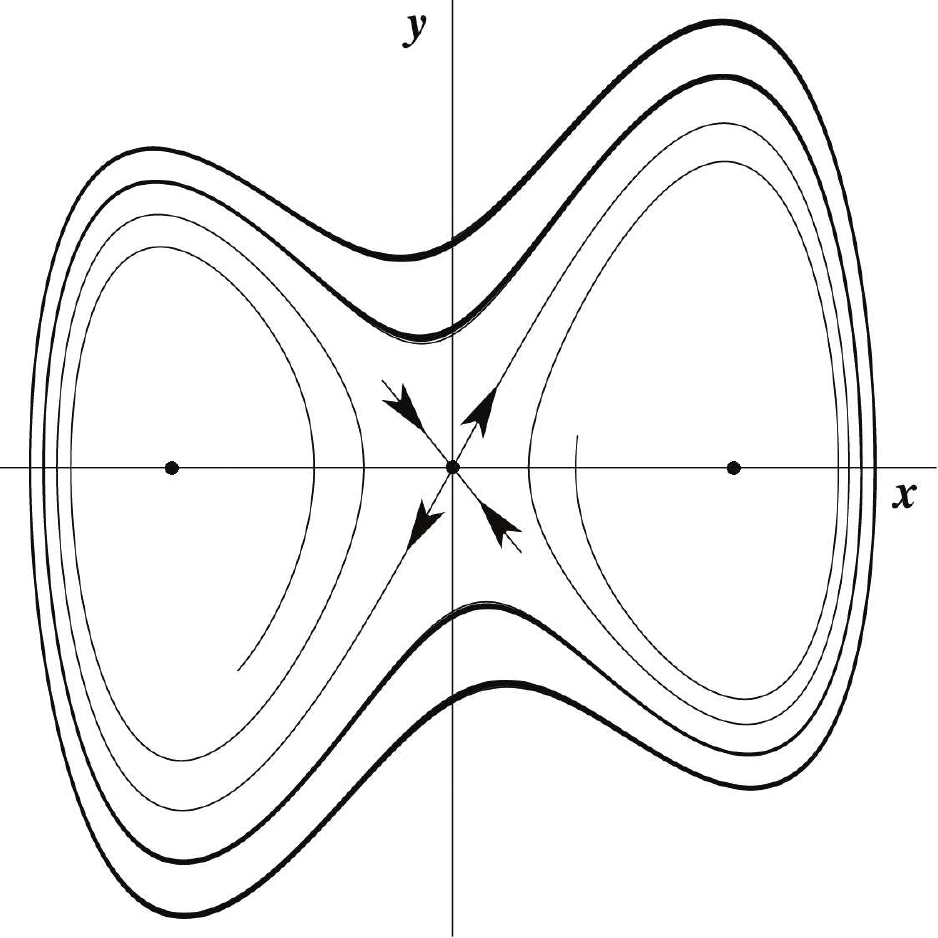}\\
\footnotesize{(g)}& \footnotesize{(h)} & \footnotesize{(i)}\\
\includegraphics[scale=0.35]{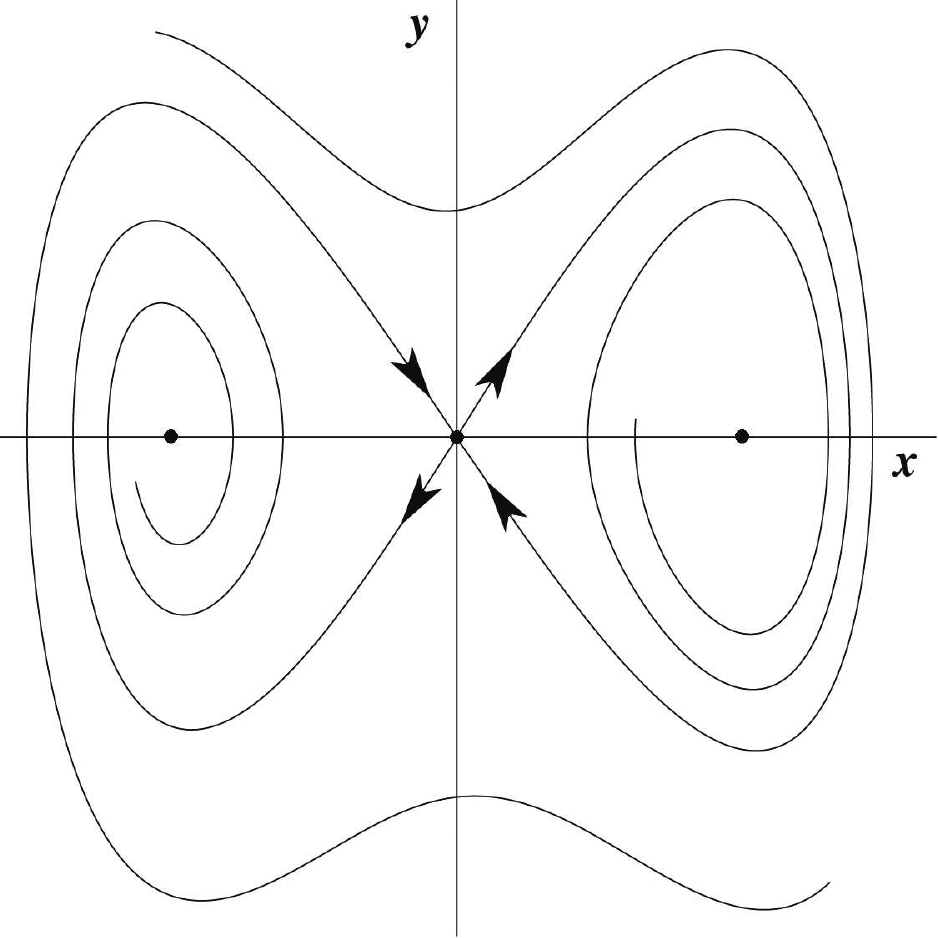}&
\includegraphics[scale=0.35]{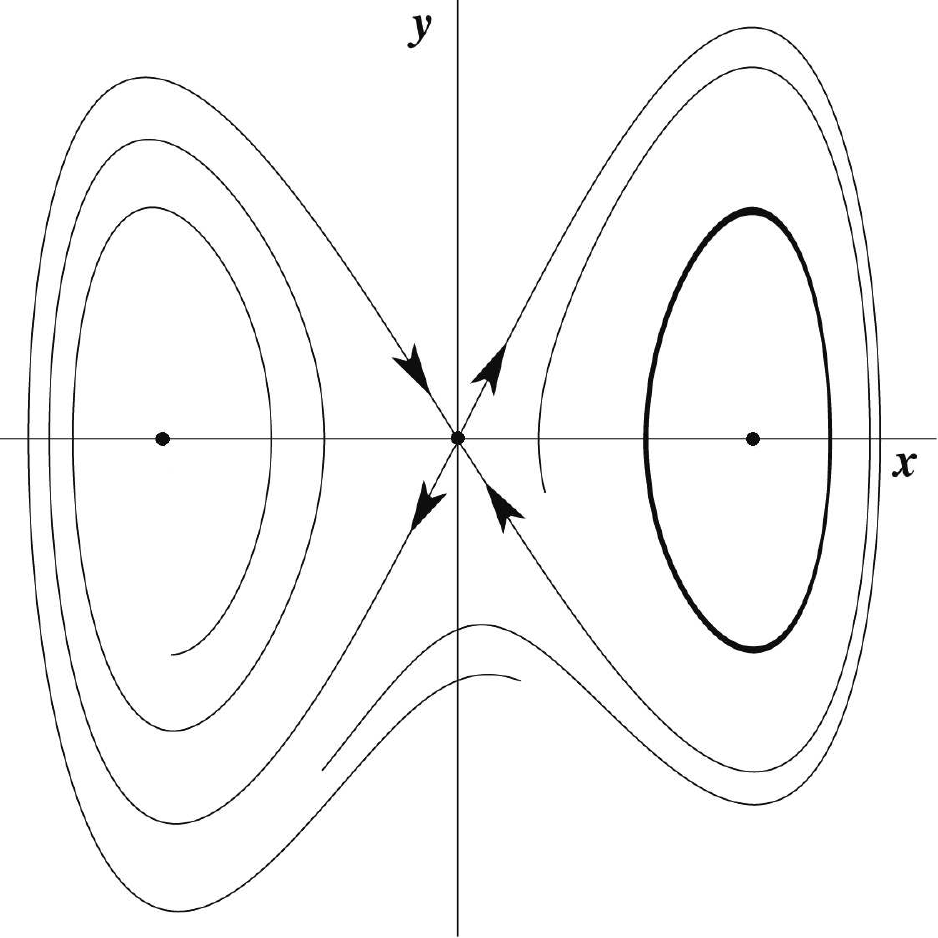}&
\includegraphics[scale=0.35]{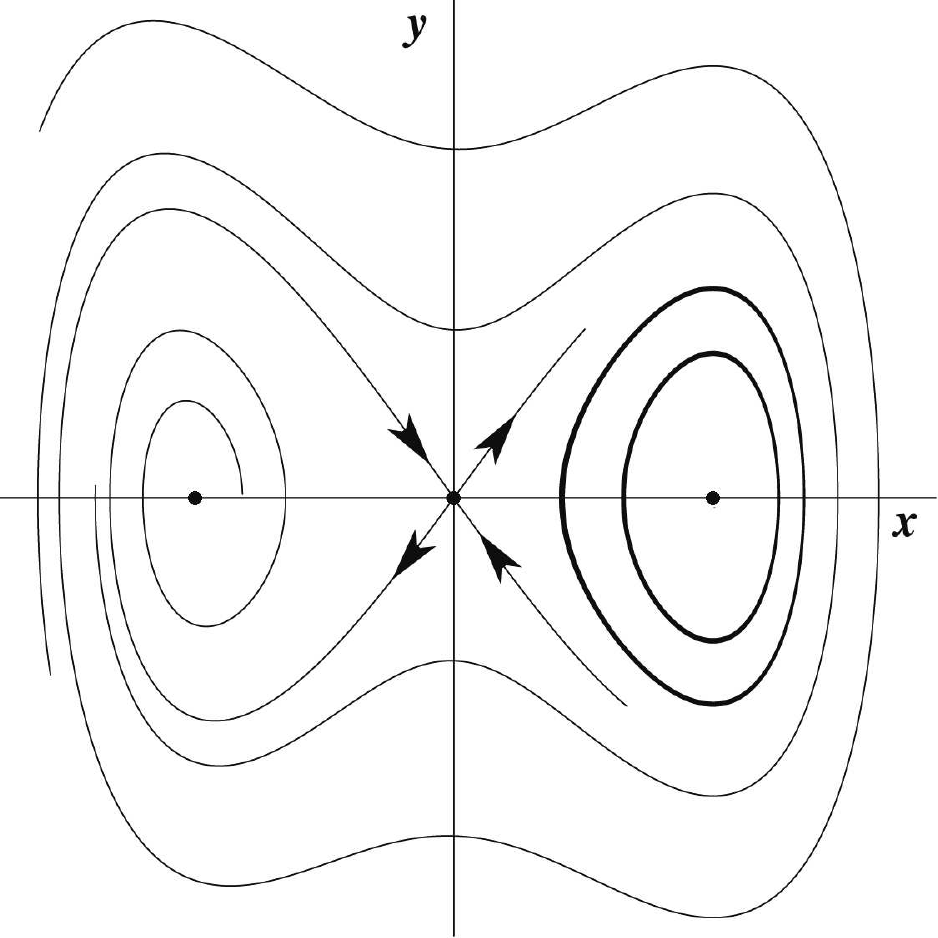}\\
\footnotesize{(j)}& \footnotesize{(k)} & \footnotesize{(l)}\\
\includegraphics[scale=0.35]{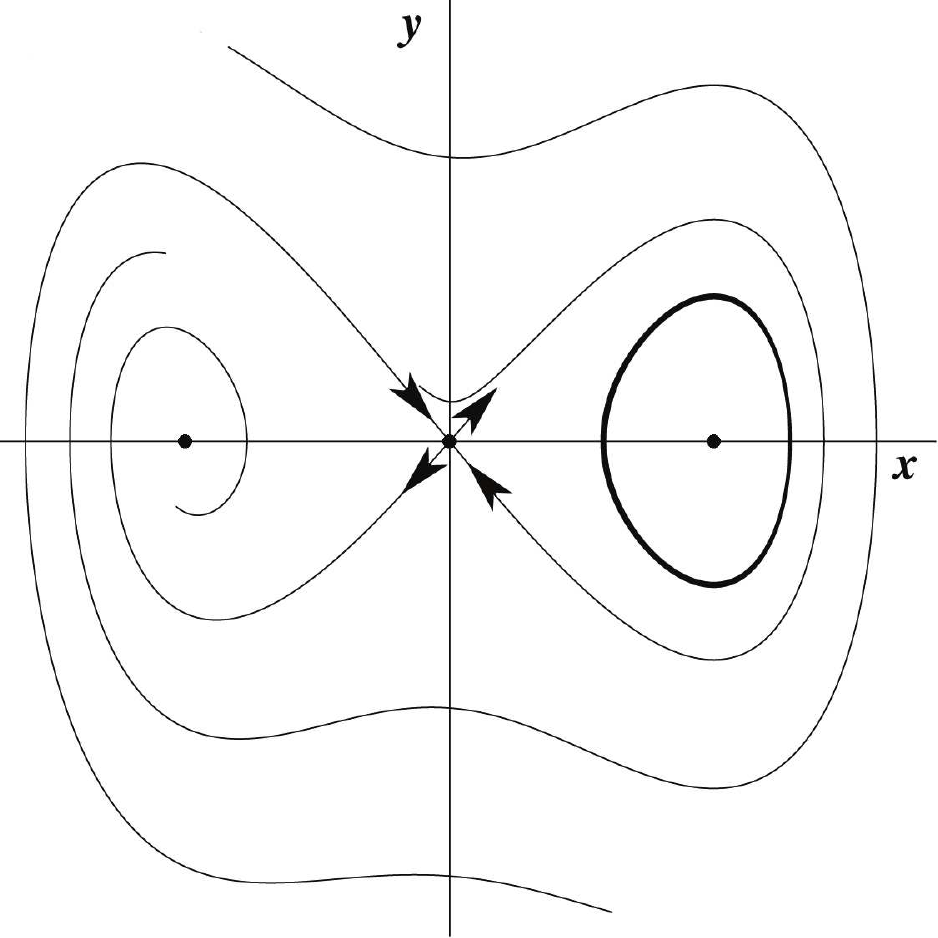}\\
\footnotesize{(m)}\\
\end{tabular}
\end{center}
\caption{Basic phase portraits of perturbed autonomous equation in
domains $D_1$ (a);  $D_2$ (b); $D_3$ (c); $D_4$ (d);  $D_5$ (e);  $D_6$ (f);  $D_7$ (g);  $D_8$ (h);  $D_9$ (i); $D_{10}$ (j);  $D_{11}$ (k);  $D_{12}$ (l); $D_{13}$ (m).}\label{fig4}
\end{figure}

Basic phase portraits of a perturbed autonomous equation for the
parameter values from 13 domains on the  $(p_1,p_2)$ plane are
presented in Fig.~\ref{fig4}. The dots show the equilibrium states
($(0,0)$ saddle point and $(\pm 1,0)$ focus), the arrows indicate
directions of motion on the separatrices. Note also that the limit
cycle near ``figure-eight'' is unstable. The other phase portraits for
symmetric domains may be obtained by rotation of angle $\pi$. The
simplest phase portraits are obtained for the dissipation domain
$D_{10}$.

\section{Analysis of  resonance zones topology}
In the domains filled with closed phase curves of the unperturbed
equation ($\varepsilon=0$) and separated from the unperturbed
separatrices, we will pass in Eq.~(\ref{eq1})  to the ``$I$
action--$\theta $ angle '' variables by the following formulas
\begin{equation}\label{eq2}
\begin{split}
 &I(h)={1\over 2\pi}\oint y(x,h)dx,\\
 &\theta={\partial S(x,I)\over \partial I}, S=\int_{x_0}^x y(x,h(I))dx,
\end{split}
\end{equation}
where $S(x,I)$ is the generating function of this canonical
transformation. The resulting system will be written in the form
\begin{equation}\label{eq2}
\begin{split}
\dot{I} &= \varepsilon
[(p_1+p_2x-x^2)y+p_3\sin{\varphi}]x^{\prime}_{\theta} \equiv \varepsilon F_1(I,\theta , \varphi ),\\
\dot{\theta} &=\omega (I) + \varepsilon
[(p_1+p_2x-x^2)y+p_3\sin{\varphi}]x^{\prime}_{I}\equiv \omega (I)+
\varepsilon F_2(I,\theta , \varphi ) , \\
\dot{\varphi} &=p_4 ,
\end{split}
\end{equation}
where  $\omega$ is the frequency of self-excited oscillations.
Consider the resonance case when
\begin{equation}\label{eq3}
   \omega (I_{pq})=(q/p)p_4 ,
\end{equation}
where $p,q$ are coprime integer numbers. The level $I=I_{pq}$
(closed phase curve  $H(x,y)=h_{pq}$ of the unperturbed system)
will be referred to as the  resonance level. The neighborhood
 $U_{\sqrt{\varepsilon}}=\{(I, \theta ): I_{pq}-C\sqrt{\varepsilon}<I<I_{pq}+C\sqrt{\varepsilon}, \ 0\leq \theta <2\pi, \  C=const>0 \}$ will be called the resonance zone.

By the substitution
\begin{equation}\label{eq4}
  \theta =  \psi +(q/p)\varphi , \quad I=I_{pq} + \mu \eta, \quad
  \mu =\sqrt{\varepsilon},
\end{equation}
in (\ref{eq2}), by averaging the obtained system over the fast
variable  $\varphi $ and neglecting the terms $O(\mu ^3)$, we
obtain the system \cite{Mor1998}
\begin{equation}\label{eq5}
    \begin{aligned}
&\dot{u } = \mu A_{0}(v,I_{pq}) + \mu ^{2}P_{0}(v,I_{pq})u, \\
& \dot{v} = \mu bu  + \mu ^{2}(b_{1}u^{2}+Q_{0}(v,I_{pq})),
\end{aligned}
\end{equation}
where $u=\eta +O(\mu ), v=\psi +O(\mu ^2)$, $b= d\omega
(I_{pq})/dI, b_1=d^2\omega (I_{pq})/2dI^2$,
\begin{equation}\label{eq6}
A_{0}(v,I_{pq})={1\over 2\pi p}\int^{2\pi
p}_{0}F_{1}(I_{pq},v+q\varphi /p,\varphi )d\varphi ,
\end{equation}
\begin{equation}\label{eq7}
P_{0}(v,I_{pq})={1\over 2\pi p}\int^{2\pi p}_{0}[\partial
F_{1}(I_{pq},v+q\varphi /p, \varphi)/\partial I]d\varphi ,
\end{equation}
\begin{equation}\label{eq8}
Q_{0}(v,I_{pq})={1\over 2\pi p}\int^{2\pi
p}_{0}F_{2}(I_{pq},v+q\varphi /p,\varphi )d\varphi .
\end{equation}

The substitution $u \to u - \mu Q_0(v, I_{pq})/b$  and the
transition to ``slow time''  $\tau =\mu t$  reduces Eqs.
~(\ref{eq5}) to a pendulum equation \cite{Mor1998}
\begin{equation}\label{eq9}
{d^{2}v\over d\tau { } ^{2}} - bA_{0}(v,I_{pq}) = \mu \sigma
(v,I_{pq}){dv\over d\tau },
\end{equation}
where
\begin{equation}\label{eq10}
    \sigma (v,I_{pq})={1\over 2\pi p}\int^{2\pi
p}_{0}(p_1 +p_2x-x^2)_{\left|{\begin{array}{c}x=x(I_{pq},v+q\varphi /p)\\
y=y(I_{pq},v+q\varphi /p)\end{array} }\right.}d\varphi.
\end{equation}
Apparently, $\sigma =const $.

The  topology of individual resonance zones may be found  from
Eq.~(\ref{eq9}) to an accuracy of terms of order $\mu^2$.

In calculations of the function $A_{0}(v,I_{pq})$ and quantities
$b$ and $\sigma$ we distinguish the following cases.

 \bfseries{Case 1} \normalfont: $(x,y) \in
G_1^{\pm}=\{(x,y):y^2/2-x^2/2+x^4/4=h, h \in (-0.25,0)\}$;

 \bfseries{Case 2} \normalfont: $(x,y) \in
G_2=\{(x,y):y^2/2-x^2/2+x^4/4=h, h>0\}$.

The unperturbed solution in (\ref{eq2}) is different in each case
\cite{KM}.

We represent the function $A_{0}(v,I_{pq})$ in the form
 $A_{0j}(v,I_{pq})=\widetilde{A}_{0j}(v,I_{pq})+B_j(I_{pq})$ and designate $b=b_j, \sigma = \sigma_j, j=1,2 $, where $B_1, B_2$ are the Poincar\'{e}--Pontryagin generating functions (\ref{eq6-1}) and (\ref{eq7-1}), respectively.

Following \cite{Mor1998}, we refer to the resonance level
$I=I_{pq}$ as splittable if the equation $A_{0j}(v;I_{pq})=0$ has
simple roots. The nonsplittable resonance level $I=I _ {pq} $ for
which $|A_{0j}(v;I_{pq})|>0$ is called passable. The splittable
resonance level $I=I _ {pq}$ is called partially passable, if
$B_j(I_{pq})\neq 0$ and impassable, if  $B_j(I_{pq})=0$.

Note that the behavior of solutions of the initial equation
(\ref{eq1}) in the neighborhood of passable, partially passable
and impassable resonance levels is defined by the theorems from
\cite{Mor1998}.

\subsection{Case 1}
Using the unperturbed solutions at the resonance level and formulas (\ref{eq6}), (\ref{eq10}), we find for $q=1$
\begin{equation}\label{eq15-1}
{d^{2}v\over d\tau ^{2}} - b_1(p_3A_1\cos{pv}+B_1) = \mu \sigma_1
{dv\over d\tau },
\end{equation}
where
\begin{equation}\label{eq13}
b_1 =\frac{\pi}{2}\frac{(2-\rho)^{3/2}[2(1-\rho){\bf
K}(\rho)-(2-\rho){\bf E}(\rho)]} {\rho^2(1-\rho){\bf K}^2(\rho)},
\end{equation}
\begin{equation}\label{eq14}
\sigma_1 = p_1 -\frac{2}{(2-\rho){\bf K}(\rho)}{\bf E}(\rho),
\end{equation}
\begin{equation}\label{eq15}
A_1=-\sqrt{2}p_4\frac{a^p}{1+a^{2p}}, a=exp\left(-\pi\frac{{\bf
K}(\sqrt{1-\rho})}{{\bf K}(\rho)} \right).
\end{equation}

For  $q>1$ we obtain the equation
\begin{equation}\label{eq15-2}
{d^{2}v\over d\tau ^{2}}  - b_1B_1 = \mu \sigma_1{dv\over d\tau }.
\end{equation}

\begin{figure}[htb]
\begin{center}
\begin{tabular}{ccc}
\includegraphics[scale=0.45]{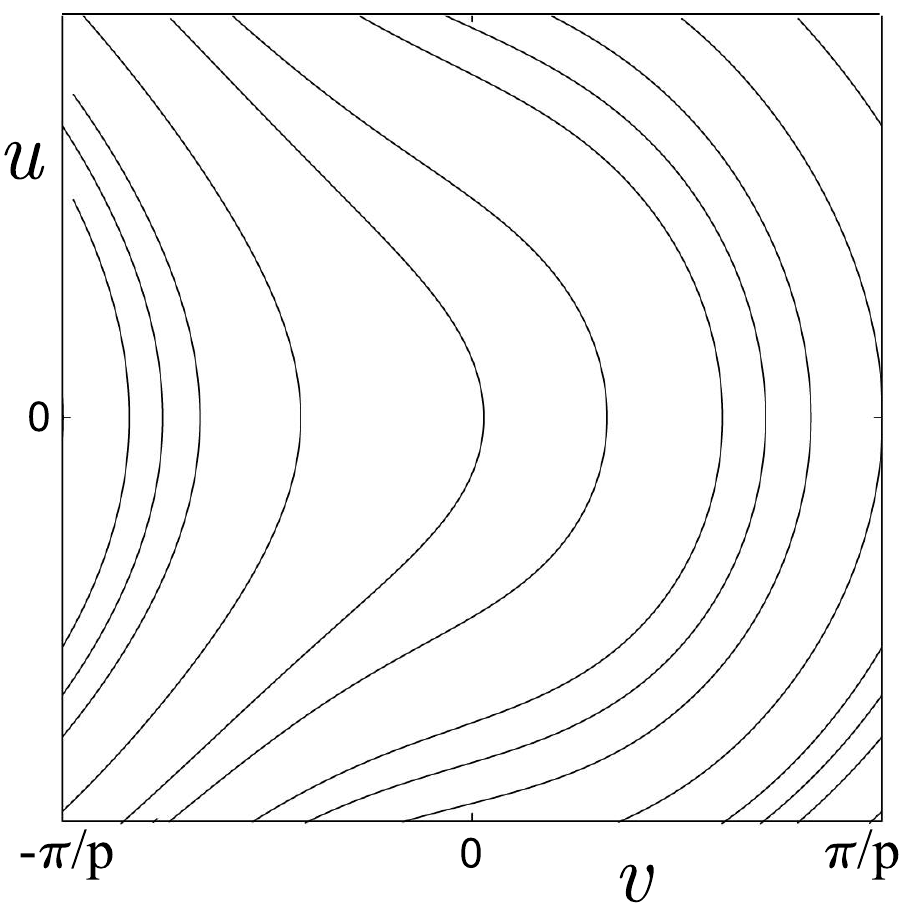}&
\includegraphics[scale=0.45]{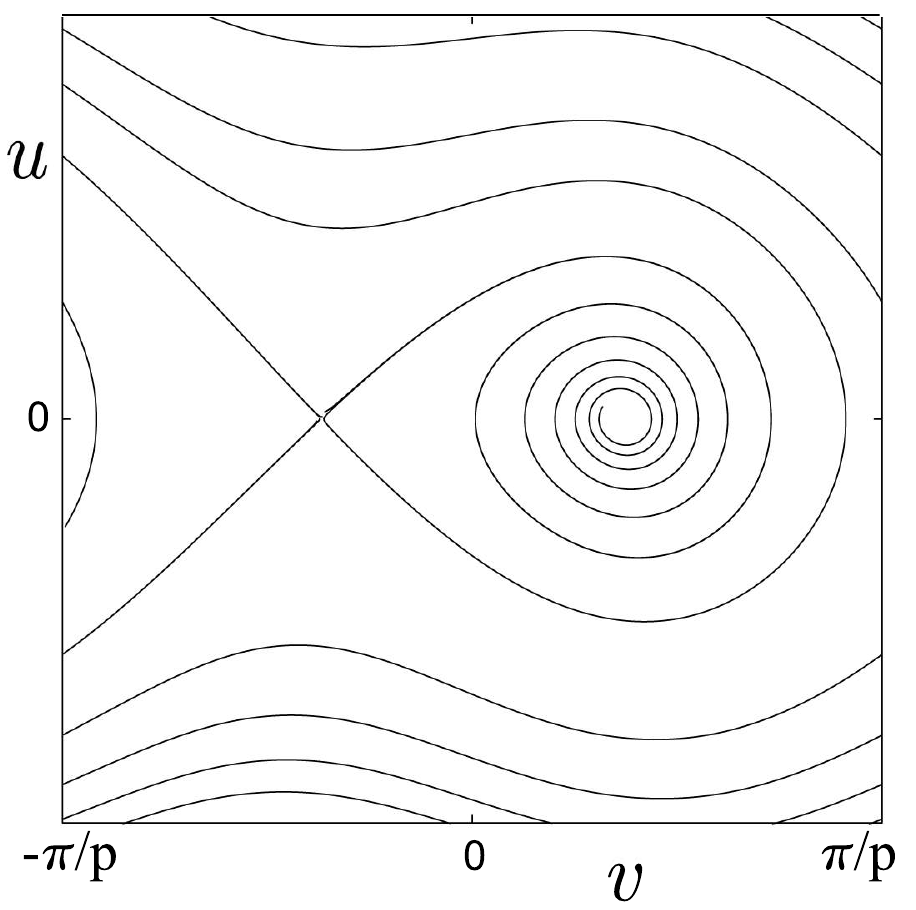}&
\includegraphics[scale=0.45]{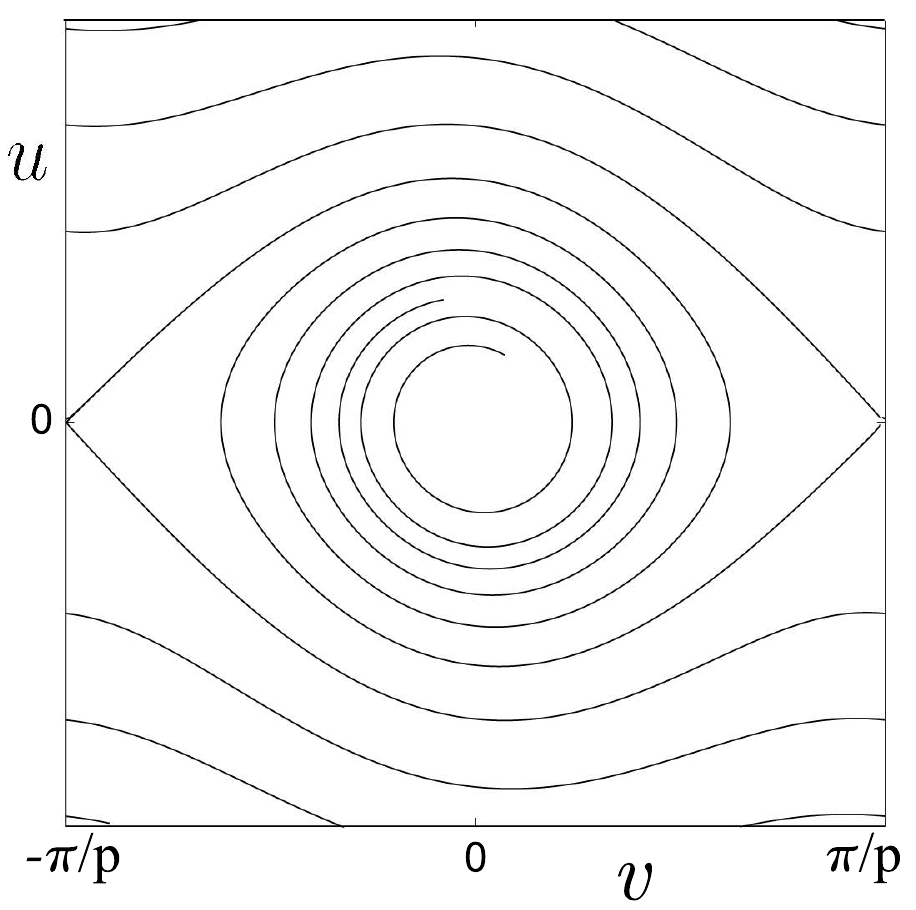}\\
\footnotesize{(a)}& \footnotesize{(b)}& \footnotesize{(c)}
\end{tabular}
\end{center}
\caption{Phase portraits for Eq.~(\ref{eq15-1}).} \label{fig5}
\end{figure}

Thus, for $q=1$ the topology of  resonance zones is described by
Eq.~(\ref{eq15-1}). The phase portraits of this equation are well
known \cite{Mor1998} (Fig. \ref{fig5}). At $B_1(I_{pq})=0$, we have {\it by
definition} an impassable resonance (Fig.~\ref{fig5}(c)). In this
case, the resonance level $I=I_{pq}$ coincides with the level
$I=I_0$ in the neighborhood of which the autonomous equation has a
limit cycle. A partially passable resonance is presented in
Fig.~\ref{fig5}(b), and a passable resonance in
Fig.~\ref{fig5}(a). According to (\ref{eq15-2}), at $q>1$ and
$B_1(I_{pq})\neq 0$ we have a passable resonance.

Consider briefly the bifurcations of the transition from the
impassable to the partially passable resonance. Let us set in
Eq.~(\ref{eq15-1}) $B_1=\mu \gamma$, where $\gamma$ defines the
deviation of the resonance level $I=I_{pq}$ from $I=I_0$. We
denote by $\gamma ^{\pm}$ the bifurcation values of  $\gamma$ at
which Eq.~(\ref{eq15-1}) has, respectively,  the upper or lower
loop enclosing a phase cylinder. As $\gamma$ deviates from the
bifurcation value, the loop gives birth to a limit cycle enclosing
a phase cylinder, and the resonance level becomes partially
passable.  The limit cycle corresponds to the two-dimensional
torus in the initial equation (for the Poincar\'{e} map it is a
closed invariant curve shown in Fig.~\ref{fig6}). More details about these
bifurcations can be found in \cite{Mor1998}.

The frequency $\omega(I)$ of the self-excited oscillations in
domains $G_1^{\pm}$ meets the condition $\omega(I)\in(0,\sqrt{2})$
and is a monotonic function. Then from the resonance condition
(\ref{eq3}) follows $p>p_4/\sqrt{2}$. Therefore, only the
resonance levels $H(x,y)=y^2/2-x^2/2+x^4/4=h_{p1}$,  for which
$p>p_4/\sqrt{2}$, are split. Note that the resonance levels with
larger values of $p$ are closer to the unperturbed separatrix.

\subsection{Case 2}
Analogously to Case 1, we find the equation
\begin{equation}\label{eq15-3}
{d^{2}v\over d\tau { } ^{2}} - b_2(p_3A_2\cos{pv}+B_2) = \mu
\sigma_2{dv\over d\tau },
\end{equation}
that defines the topology of the resonance zones    at odd $p$ and
$q=1$. Otherwise, the resonance zones  topology  is defined by the
equation
\begin{equation}\label{eq15-4}
{d^{2}v\over d\tau { } ^{2}}  - b_2B_2 = \mu \sigma_2{dv\over
d\tau }.
\end{equation}
In (\ref{eq15-3}) and (\ref{eq15-4}) we have
\begin{equation}\label{eq22}
b_2 =\frac{\pi}{4}\frac{(2\rho-1)^{3/2}[(1-\rho){\bf
K}(\rho)+(2\rho-1){\bf E}(\rho)]} {\rho(1-\rho){\bf K}^2(\rho)},
\end{equation}

\begin{equation}\label{eq23}
\sigma_2 = p_1 -\frac{2}{(2\rho-1){\bf K}(\rho)}({\bf
E}(\rho)+(\rho-1){\bf K}(\rho)),
\end{equation}
\begin{equation}\label{eq24}
A_2=-2\sqrt{2}p_4\frac{a^{p/2}}{1+a^p}.
\end{equation}

For even $p$ and/or $q>1$, the resonance is passable if
$B_2(I_{pq})\neq 0$.

\begin{figure}[!h]
\begin{center}
\begin{tabular}{cc}
\includegraphics[scale=0.6]{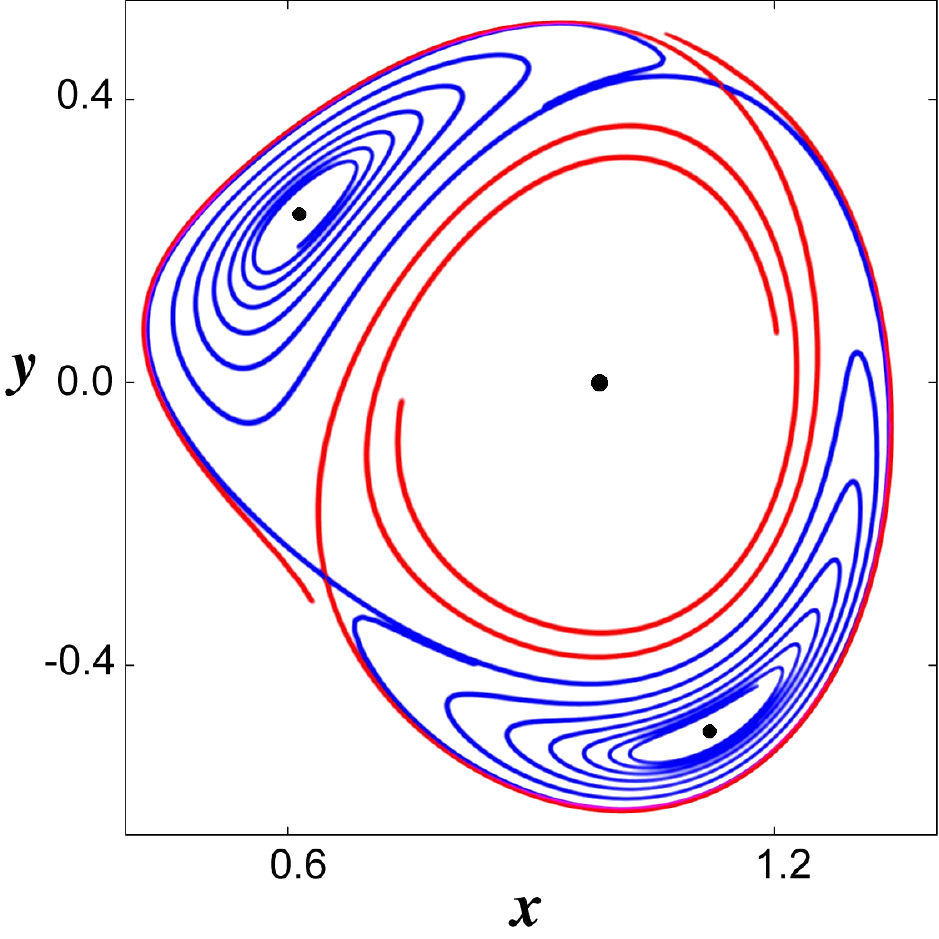}&
\includegraphics[scale=0.6]{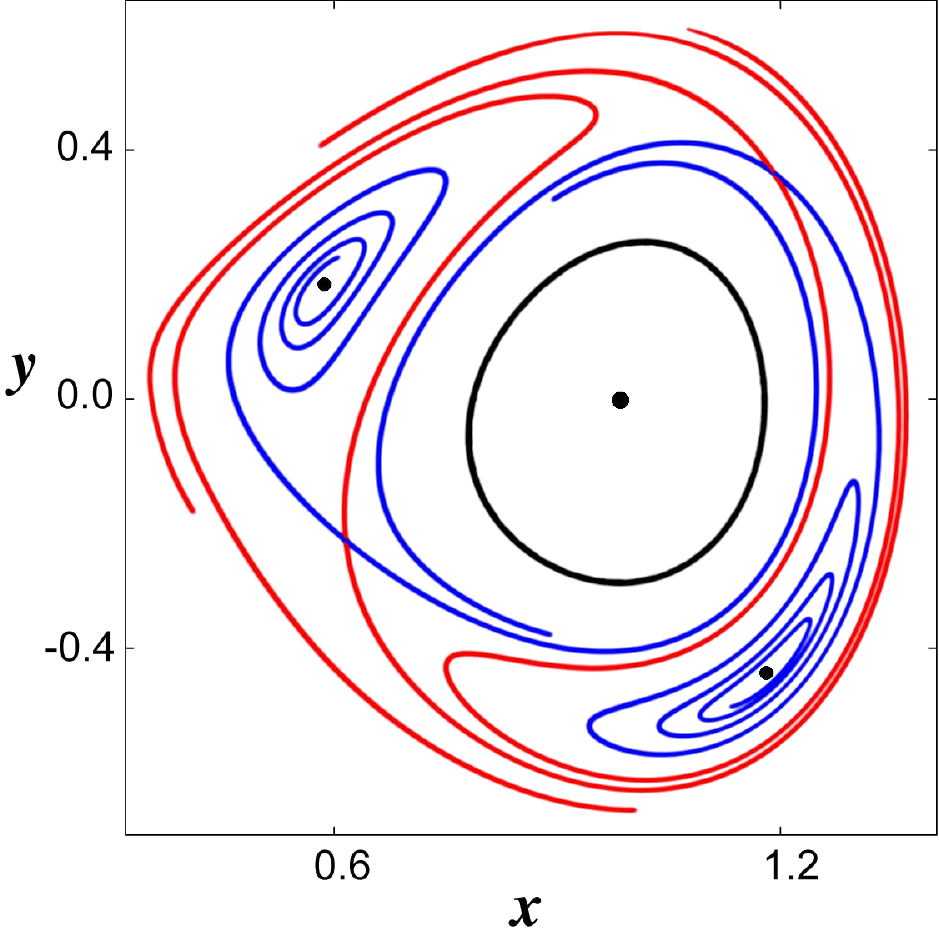} \\
\footnotesize{(a)}& \footnotesize{(b)}\\
\includegraphics[scale=0.6]{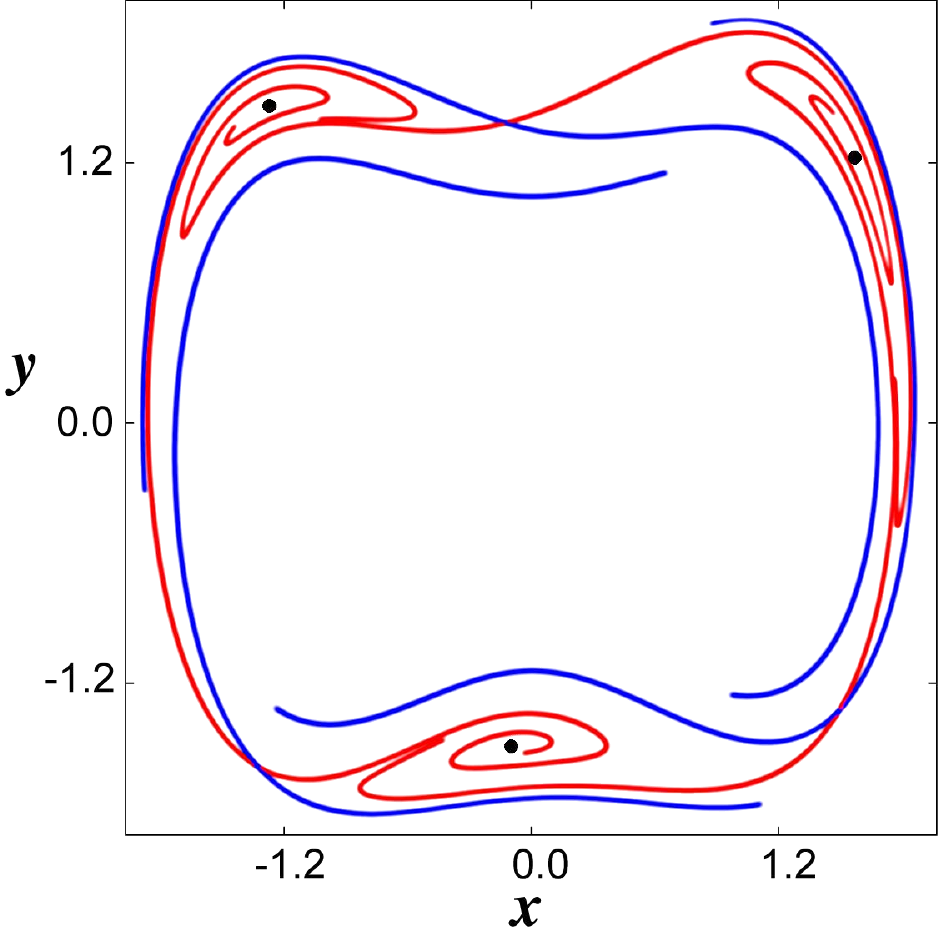}&
\includegraphics[scale=0.6]{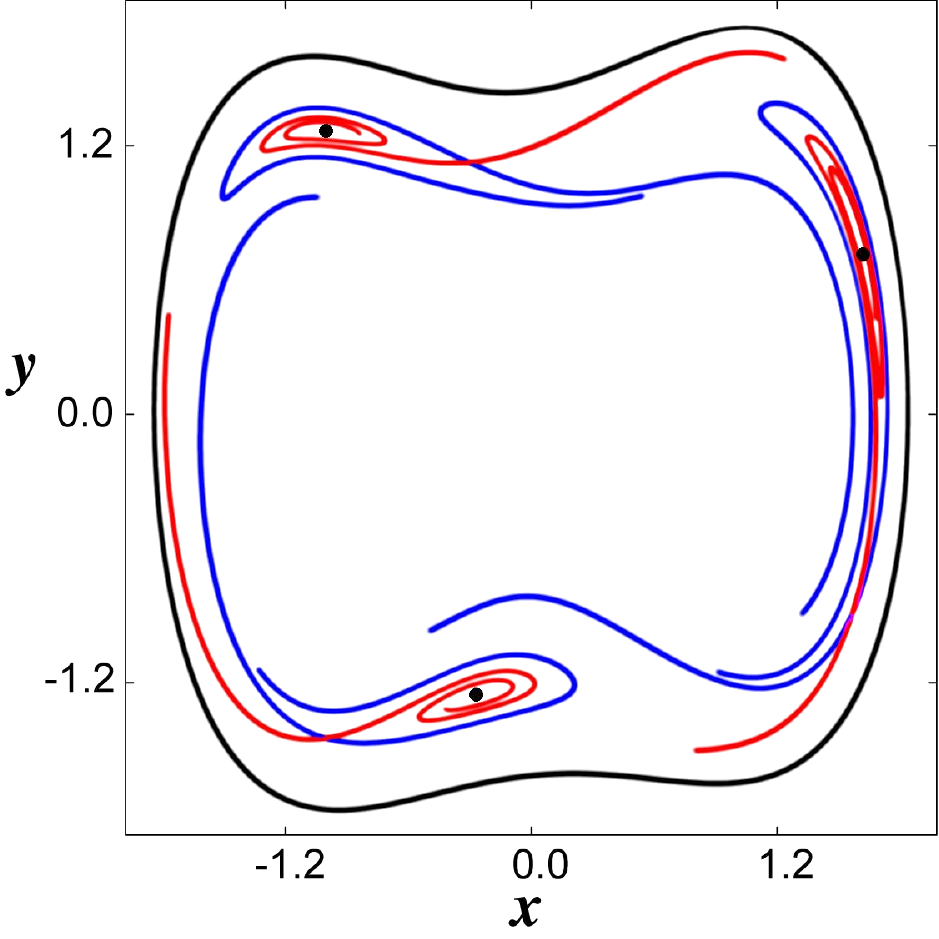}\\
 \footnotesize{(c)} & \footnotesize{(d)}\\
\end{tabular}
\end{center}
\caption{Behavior of invariant curves of the Poincar\'{e} map for
Eq.~(\ref {eq1}) at $\varepsilon=0.1$ and  $ p_1=1, p_2=-0.1,
p_3=0.5, p_4=2.5$ (a);  $ p_1=1, p_2=-0.02, p_3=0.5, p_4=2.5$ (b);
$ p_1=1, p_2=0.03, p_3=1, p_4=3.36$ (c); $ p_1=1, p_2=0.03, p_3=1,
p_4=3$ (d).}\label{fig6}
\end{figure}

The Poincar\'{e} map for Eq.~(\ref{eq1}) at different parameter
values was constructed using the WInSet software
\cite{MD}\footnote{The first version of the software was described
in \cite{MDBM}.}. It was found that at small values of
$\varepsilon$ numerical results are in a good agreement with the
theoretical study.  Figure~\ref{fig6} illustrates the structure of
the neighborhood of the splittable levels $I=I_{21}, I=I_{31}$.
The impassable resonance zones are shown in Figs.~\ref{fig6}(a) and \ref{fig6}(c), and the partially passable zones in Figs.~\ref{fig6}(b) and \ref{fig6}(d). The dots in Fig.~\ref{fig6} correspond to  points of period-2, as well as to the fixed points of the Poincar\'{e} map in domain $G_1^+$
(Figs.~\ref{fig6}(a) and \ref{fig6}(b)) and periodic points of period-3  in domain $G_2$ (Figs.~\ref{fig6}(c) and \ref{fig6}(d)). Besides, a closed invariant curve of
the Poincar\'{e} map in domain $G_1^+$ is shown in
Fig.~\ref{fig6}(b) and in domain $G_2$ in Fig.~\ref{fig6}(d). The
stable separatrices are plotted by the blue curves, the unstable
separatrices by the red ones.

Note that the fixed and periodic points in resonance zones
correspond to resonance periodic solutions of period $2\pi p /p_4$
in the initial equation, and the closed invariant curves to
quasiperiodic (double-frequency) solutions (two-dimensional tori).

\section{On global behavior of solutions outside the neighborhood of ``figure-eight''}

The resonance levels corresponding to the limit cycles in the
autonomous equation ($p_3=0$ in Eq.~(\ref{eq1})) are splittable as $B_ {j}
(I_{pq}) =0$, $j=1,2$ for these levels. According to Sec.~2, the
number of such levels in each domain $G_1^{\pm}, G_2$ is no more
than two. Let us remove from these domains the neighborhoods of
such levels and designate the remaining domains without the
neighborhood of ``figure-eight'' by $V$. Using the results obtained in
\cite{Mor1998} and Eqs.~(\ref{eq15-1}), (\ref{eq15}), (\ref{eq15-3}) and (\ref{eq24}), we obtain the following theorem.
\begin{theorem}\label{th1} There are only finitely many splittable resonance levels in $V$.

\end{theorem}

It follows from this theorem that for relatively small
$\varepsilon > 0$ the neighborhoods of splittable resonance levels
do not intersect. This allows us to speak about the global
behavior of solutions in the considered cells. According to \cite
{Mor1998}, the separatrices of saddle periodic solutions lying at
different  resonance levels intersect, causing the formation of
heteroclinic structures and a complicated geometry of the
attraction domains of stable periodic solutions.

\begin{figure}[!ht]
\begin{center}
\begin{tabular}{cc}
\includegraphics[scale=0.8]{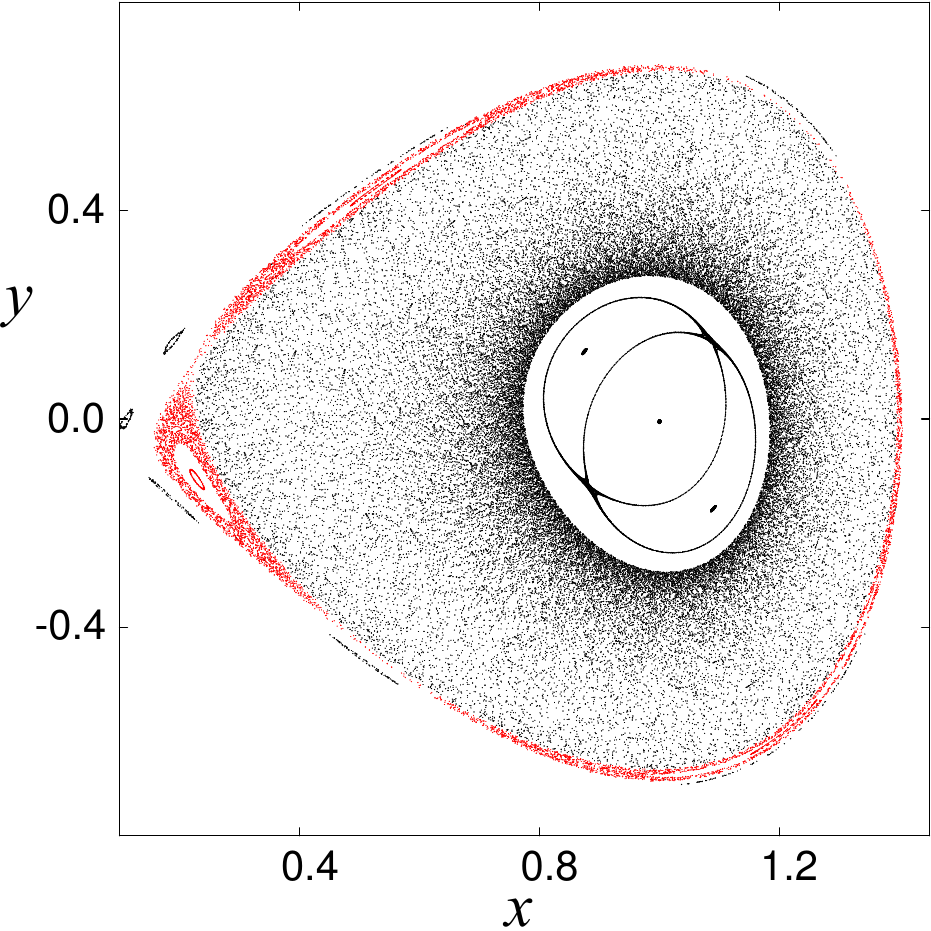}&
\includegraphics[scale=0.8]{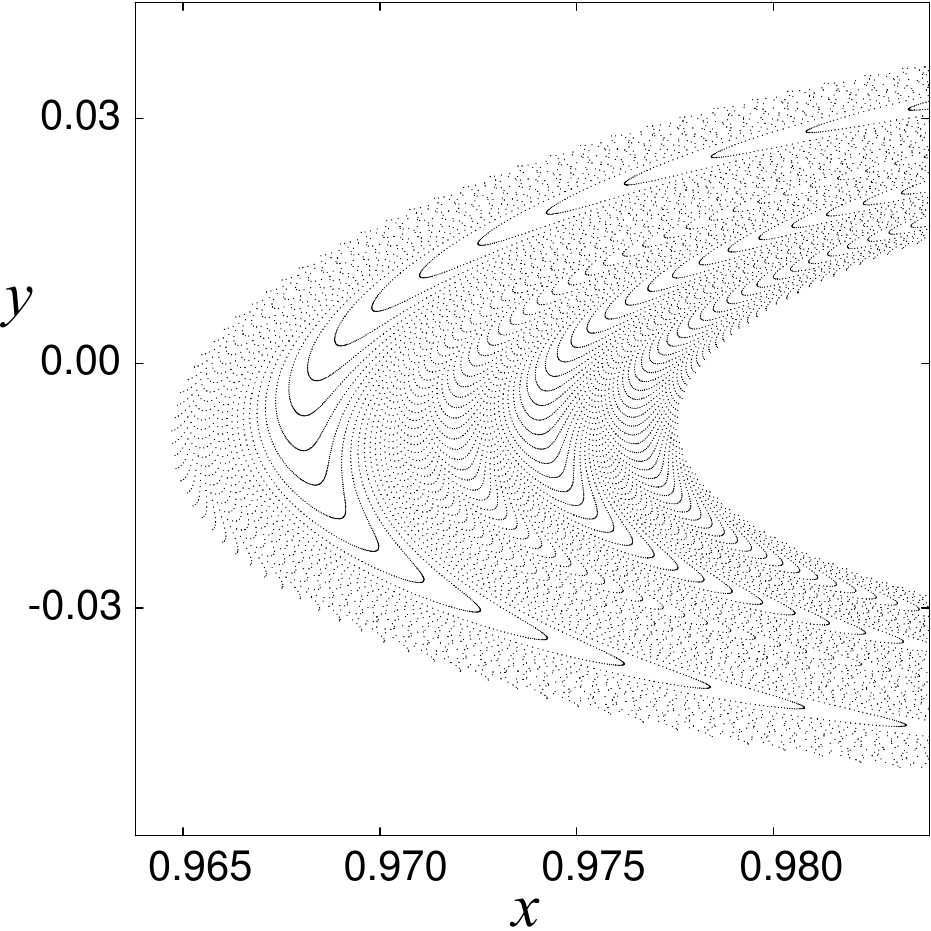}\\
\footnotesize{(a)}& \footnotesize{(b)}
\end{tabular}
\end{center}
\caption{Behavior of the trajectories of the Poincar\'{e} map for
Eq.~(\ref {eq1})  in the case of two impassable resonances with
$p=2 $ and $p=3 $ in domain $G_1 ^ {+} $ (a), and a fragment of
this domain  (b). Here, $\varepsilon=0.01$, $p_1=-0.221$,
$p_2=1.22$, $p_3=1$, $p_4=2.782$.}\label{fig7}
\end{figure}

Consider as an example of the illustration of global behavior of solutions, a
cell inside the right loop. Let the perturbed autonomous equation
have two limit cycles in this cell (domain $D_{12}$ in the
bifurcation diagram) and $ \rho =\rho ^ {(1)} $, $\rho =\rho ^
{(2)} $ be simple roots of the Poincar\'{e}--Pontryagin function
$B_1^{+}(\rho )$ in domain $G_1^{+}$.

Let us fix the parameter $p_2=1.22$ and find the values of the
parameters $p_1$, $p_4 $ at which the cycle $\rho=\rho^{(1)} $
coincides with the resonance level $I=I _ {21} $, and the cycle $
\rho =\rho ^ {(2)} $ with the level $I=I _{31}$. From the
resonance condition (\ref{eq3}) we have $p_4=2\omega (\rho ^{(1)},
p_1) =3\omega (\rho ^ {(2)},p_1)$. From this relation and the
equations defining the limit cycles $B_1^{+}(\rho^{(1)}, p_1) =0$,
$B_1 ^{+} (\rho^{(2)}, p_1) =0 $, we find $p_1\approx -0.221 $, $
\rho ^{(1)} \approx 0.45 $, $ \rho ^{(2)} \approx 0.98$. Then, at
$p_4\approx 2.782 $ we will have impassable resonance zones. There
exists in the initial equation a stable periodic solution of period
$6\pi /p_4 $ (stable periodic points of period-3 in the red zone in Fig.~\ref {fig7}(a)) and an unstable periodic solution of period $4\pi /p_4 $ (unstable periodic points of period-2  in the white zone in Fig.~\ref{fig7}(a)). All the other resonance levels between the above ones will be passable.
Between the resonance zones, trajectories fill the cell under
consideration (see Fig.~\ref {fig7}(a); actually, these
trajectories are slowly spiraling and tend to stable periodic
points of period-3 as  $t\to \infty $). A higher-order partially
passable resonance can be seen near the unperturbed separatrix
loop in Fig.~\ref {fig7}(a). A magnified fragment with a
trajectory inside the resonance zone with $p=2$ is presented in
Fig.~\ref{fig7}(b), where one can see passable resonances.
Passable resonances were also observed between the  resonance
zones with $p=2$ and $p=3$.

The behavior of the invariant curves of the Poincar\'{e} map in
the neighborhood of the impassable resonance zone with $p=2$ is
shown in more detail in Fig.~\ref{fig6}(a).

\section{Analysis of the behavior of solutions in the small neighborhood of ``figure-eight''}

The unperturbed equation $ \ddot {x} -x+x^3=0 $ has a right loop
  $ \Gamma ^{r} = \Gamma ^{r} _ {s} \bigcup \Gamma ^{r}_{u} $ of
saddle separatrix $O (0,0) $ and the left loop $\Gamma ^{l} =
\Gamma ^{l}_{s} \bigcup \Gamma ^{l}_{u}$ (Fig. \ref{fig1}).

It is known that under the action of perturbations, the
separatrices of the fixed saddle point of the Poincar\'{e} map may
intersect forming homoclinic structures of two types:
 1) $\Gamma^{r}_{s}\bigcap \Gamma^{r}_{u}\neq \oslash $ and/or
 $\Gamma^{l}_{s}\bigcap \Gamma^{l}_{u}\neq \oslash $;
 2)  $\Gamma^{r}_{u}\bigcap \Gamma^{l}_{s}\neq \oslash $
or $\Gamma^{r}_{s}\bigcap \Gamma^{l}_{u}\neq \oslash $, when
$p_2\neq 0$.

Existence of a homoclinic structure results in complicated
behavior of solutions in its neighborhood or, in other words, in a
nontrivial hyperbolic set \cite{Sh1967}. The problem of the
existence of type 1)  homoclinic structure is solved using the
Melnikov formula \cite {Melnikov} $\Delta (t_0) = \varepsilon
\Delta _1 (t_0) +O (\varepsilon^2)$, where $ \Delta (t_0)$ is the
distance between the related branches of the separatrix into which
the unperturbed separatrix splits. The substitution of  $x =\xi
+\varepsilon x_1 (t) +O({\varepsilon} ^2) $, where
\begin{equation}\label{eq30}
x_1(t)=-\frac{p_3}{1+{p_4}^2}\sin{(p_4 t)},
\end{equation}
in (\ref{eq1}) yields the following equation
\begin{equation}\label{eq31}
\ddot{\xi}-\xi+{\xi}^3=\varepsilon \left
[(p_1+p_2\xi-{\xi}^2)\dot{\xi}+\frac{3p_3}
{1+{p_4}^2}{\xi}^2\sin{(p_4 t)}\right]. \\
\end{equation}
Applying the Melnikov formula to this equation, we find
 \begin{equation}\label{eq32}
\Delta_1(t_0)=2\left(\frac{2}{3}p_1\pm \frac{\pi}{8}\sqrt{2}p_2 -
\frac{8}{15}\right)+ \frac{3\pi p_4}{2{\rm cosh}(\pi
p_4/2)}p_3\cos{(p_4 t_0)}.
\end{equation}
If $ \Delta_1 (t_0) $ is an alternating function, which holds
under the condition
\begin{equation}\label{eq33}
|p_3| \ > p_3^* = \frac{4}{3}\left|\left(\frac{2}{3}p_1 \pm
\frac{\pi}{8}\sqrt{2}p_2 - \frac{8}{15}\right)\frac{{\rm cosh}(\pi
p_4/2)}{\pi p_4}\right|,
\end{equation}
then there occurs transversal intersection of the stable and
unstable manifolds of the fixed point.

If $ \Delta_1 (t_0) $ is a constant-sign function, then the
corresponding separatrix manifolds of the saddle fixed point do
not intersect. However, if the value of $|p_3-p_3^*|$ is small
enough, then, as follows from \cite{GSh}, \cite{Mor1976},  a
nontrivial hyperbolic set exists in the neighborhood of ``figure-eight''.

Under the condition
\begin{equation}\label{eq34}
|p_3| \ = \ \frac{4}{3}\left|\left(\frac{2}{3}p_1\pm
\frac{\pi}{8}\sqrt{2}p_2 - \frac{8}{15}\right)\frac{{\rm
cosh}{(\pi p_4/2)}}{\pi p_4}\right|
\end{equation}
the corresponding separatrices of the fixed point $(0, 0)$ are
tangent to each other (to an accuracy of terms of order
$\varepsilon ^2$).

Making use of the Melnikov formula, it is easy to represent all
possible cases of relative position of the separatrices as a
result of splitting of the left or right separatrix loop. For
example, for $p_3=0$, the condition
$$\frac{2}{3}p_1+\frac{\pi}{8}\sqrt{2}p_2 - \frac{8}{15}=0$$
specifies the existence of the right separatrix loop. With allowance
for external force, the outgoing and incoming separatrices
intersect transversally, forming a homoclinic Poincar\'{e}
structure.  In this case, for the left separatrix loop we have
\begin{equation}\label{eq35}
\Delta_1(t_0)= -\frac{\pi}{2}\sqrt{2}p_2+\frac{3\pi p_4}{2{\rm
cosh}{(\pi p_4/2)}}p_3\cos{(p_4 t_0)}.
\end{equation}

\begin{figure}[htb]
\begin{center}
 \begin{tabular}{ccc}
\includegraphics[scale=0.4]{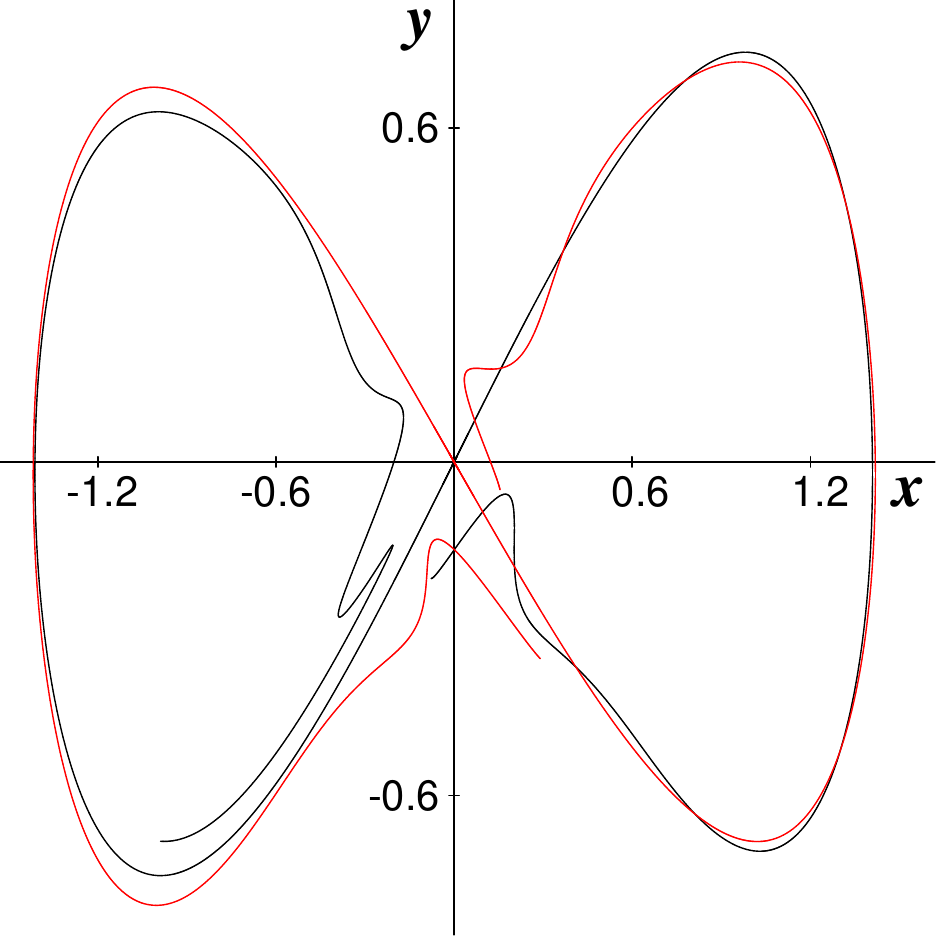}&
\includegraphics[scale=0.4]{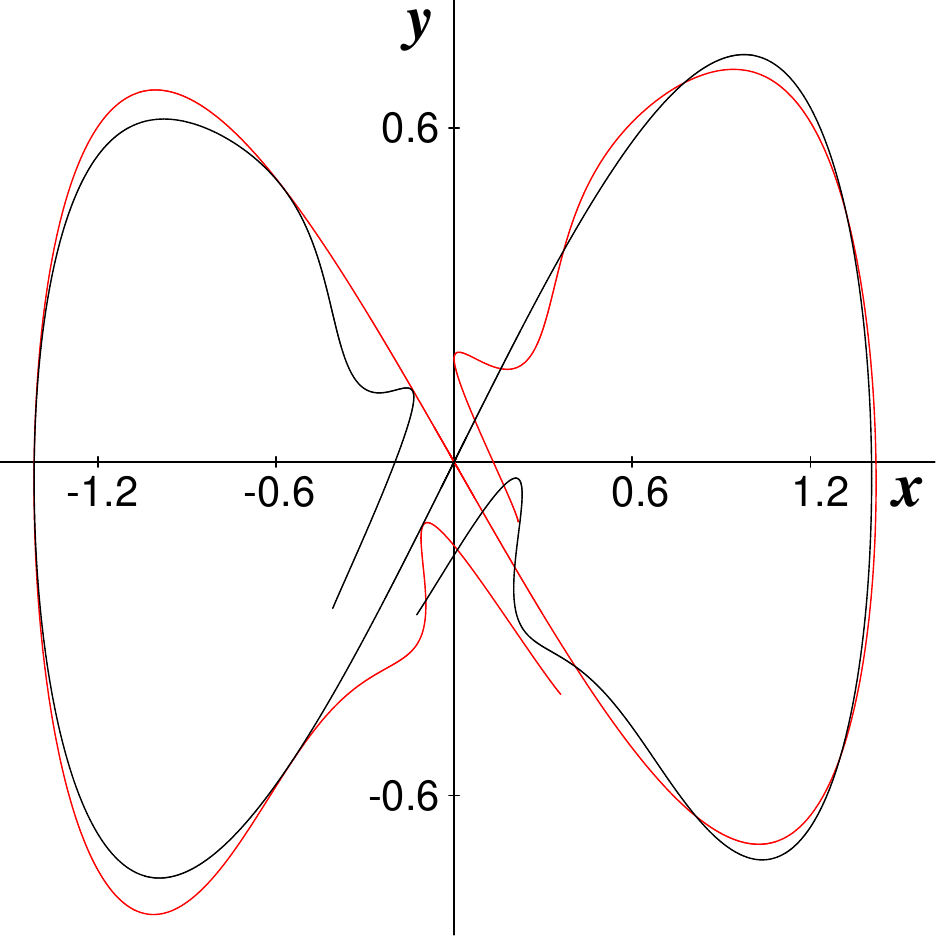}&
\includegraphics[scale=0.4]{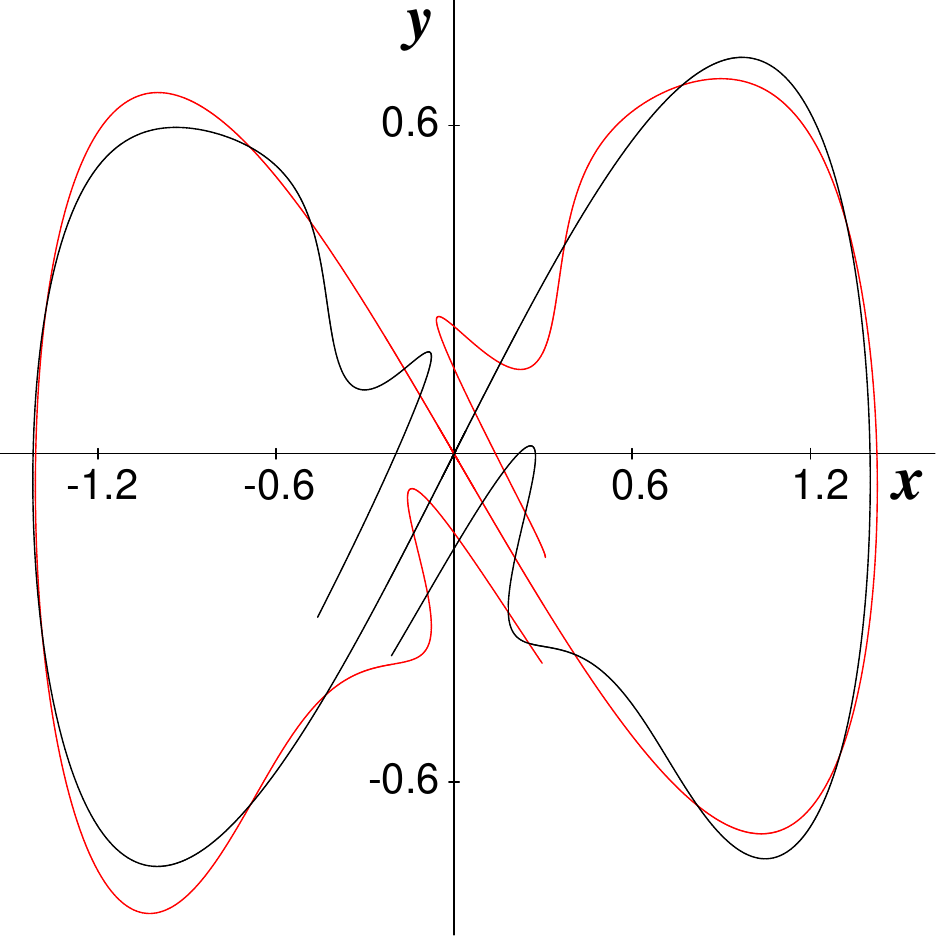}\\
\footnotesize{(a)} & \footnotesize{(b)} & \footnotesize{(c)}\\
 \end{tabular}
\end{center}
\caption{Behavior of separatrices of the fixed point $(0,0)$ for
Eq.~(\ref {eq31}) on the $ (\xi = x, \dot {\xi} = y) $ plane at  $\varepsilon=0.3$, $p_1=0.7551195621$,  $p_2=0.053875454$,
$p_4=4$ and $p_3=1.13$ (a),  $p_3=1.7$ (b),  $p_3=2.83$ (c).}
\label{fig8}
\end{figure}

\begin{figure}[htb]
\begin{center}
 \begin{tabular}{cc}
\includegraphics[scale=0.45]{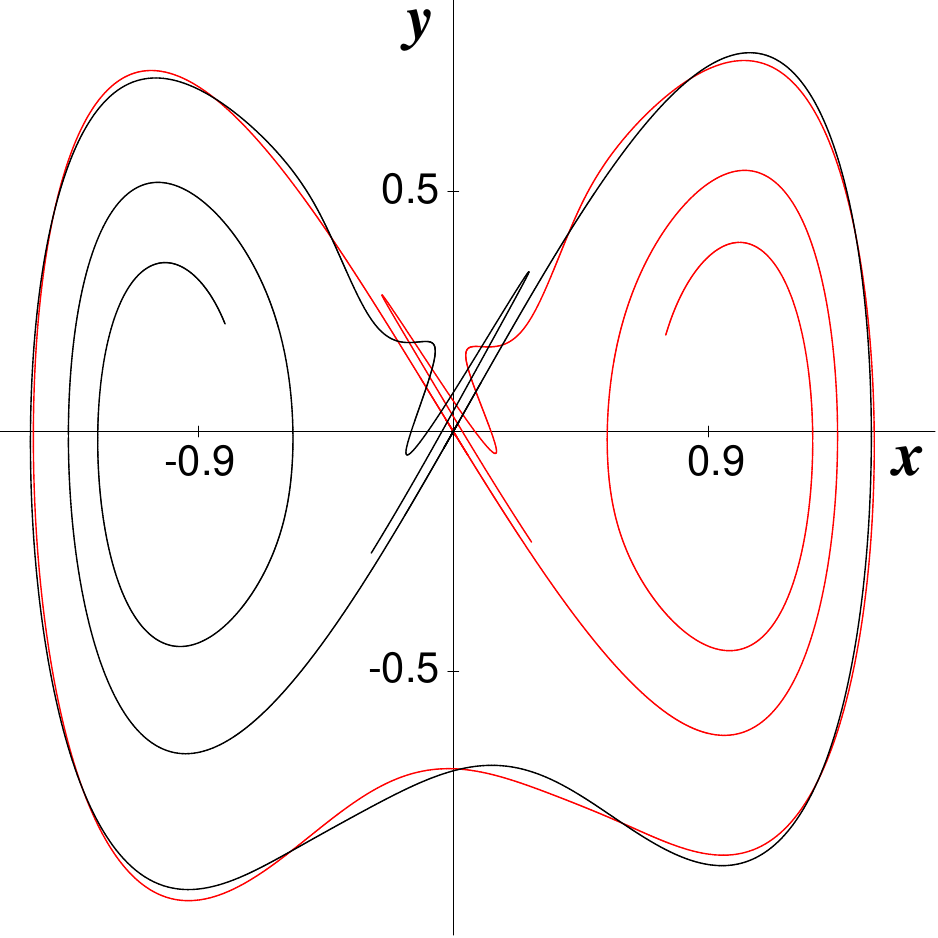}&
\includegraphics[scale=0.45]{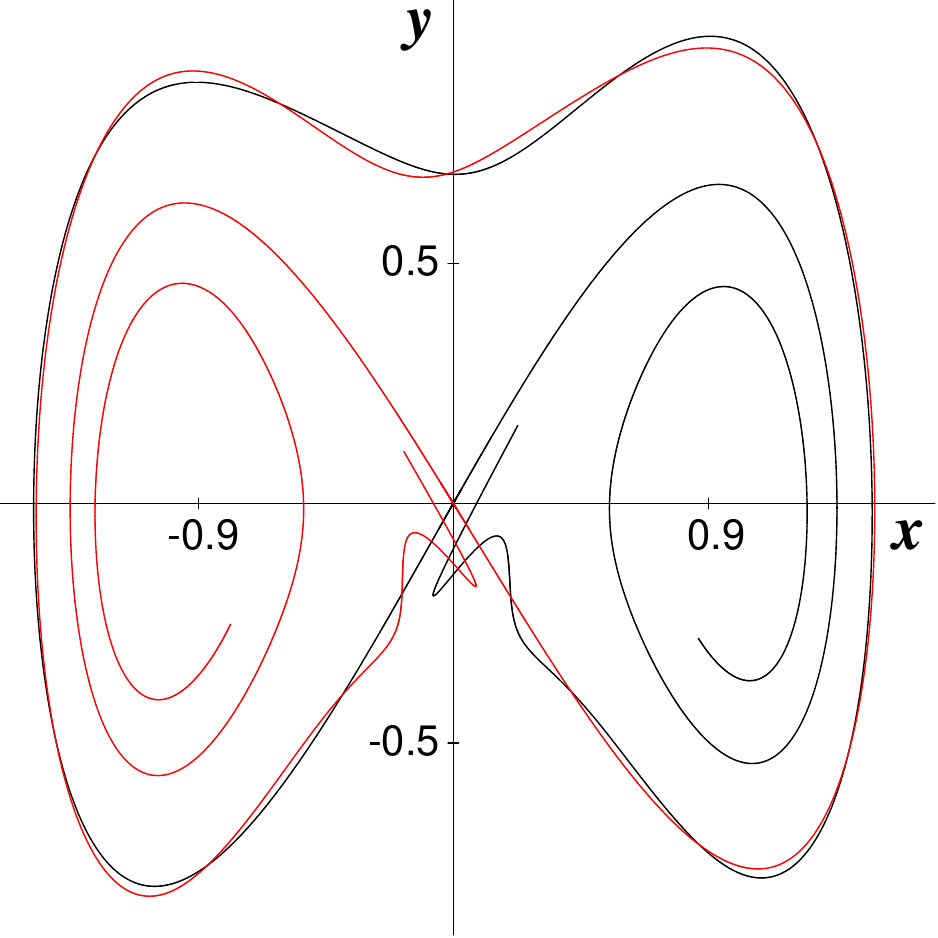}\\
\footnotesize{ (a) } & \footnotesize{ (b) }\\
 \end{tabular}
\end{center}
\caption{Behavior of separatrices of the fixed point $(0,0)$ for
Eq.~(\ref{eq31}) at  $\varepsilon=0.1$, $p_1=0.78549$, $p_3=1.02$,
$p_4=4$ and (a) $p_2=1.6$, (b) $p_2=-1.6$.}\label{fig9}
\end{figure}

\begin{figure}[htb]
  \begin{center}
    \begin{tabular}{cc}
    \includegraphics[scale=0.4]{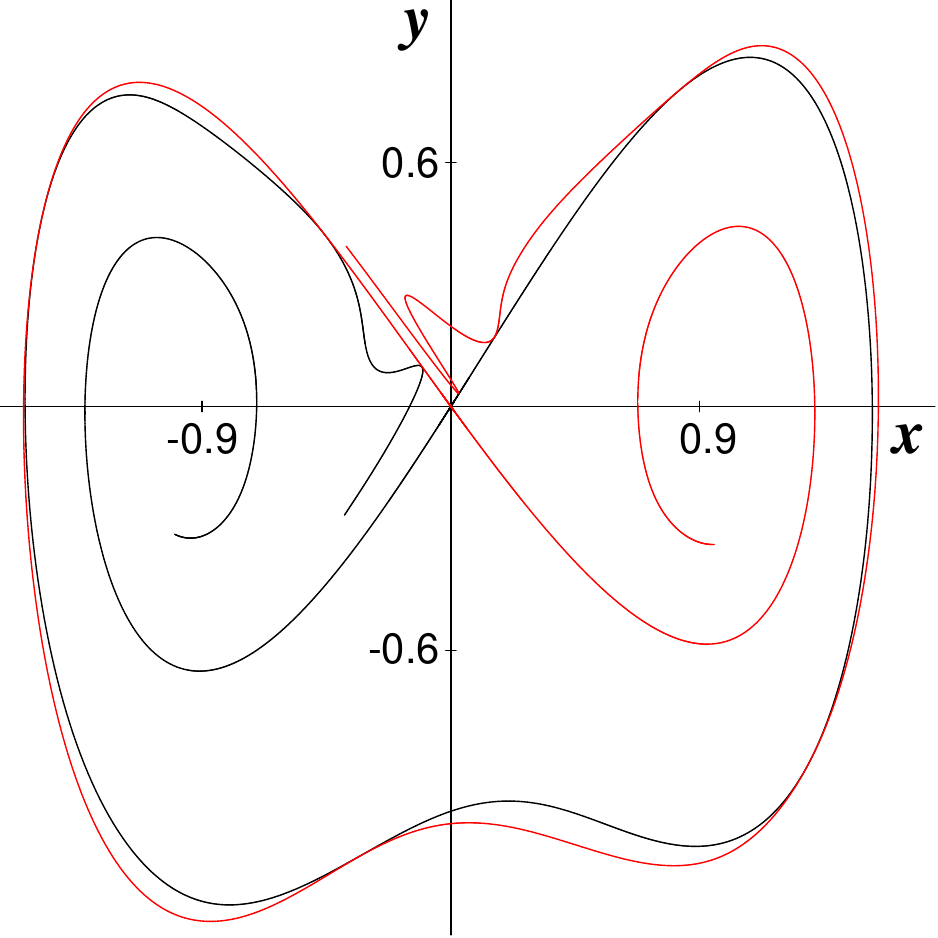} &
      \includegraphics[scale=0.4]{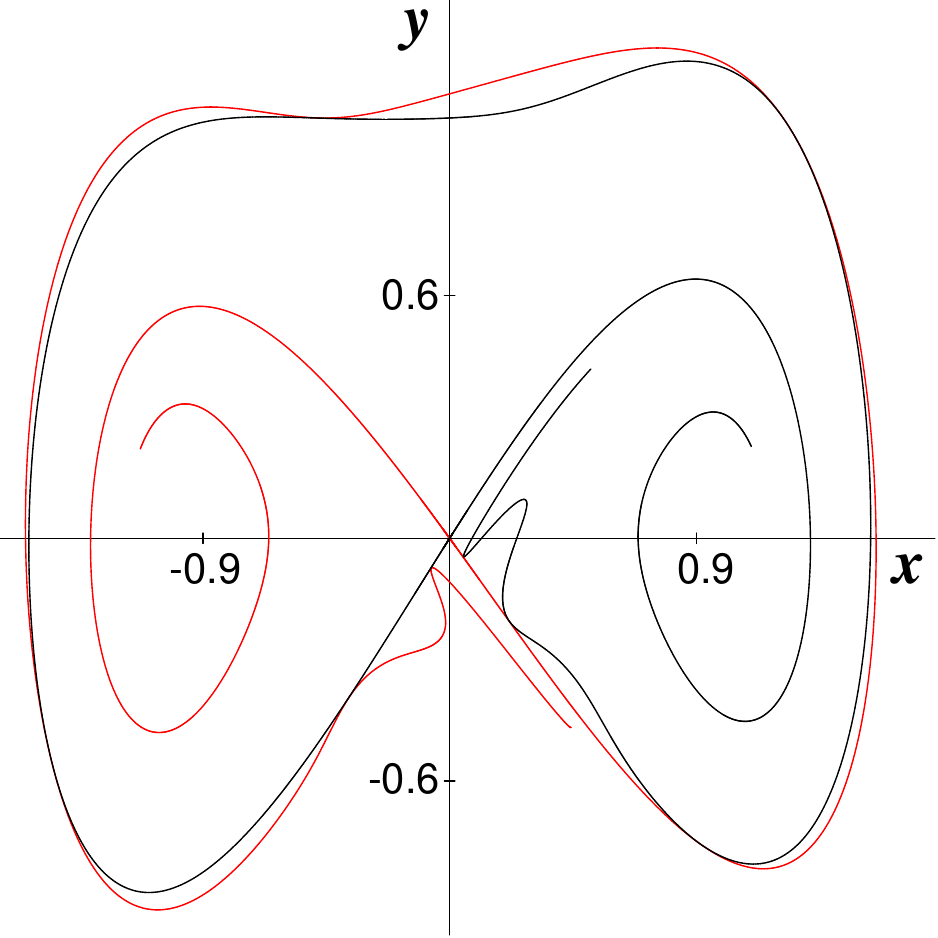}\\
       \footnotesize{(a)}& \footnotesize{(b)}\\
      \includegraphics[scale=0.4]{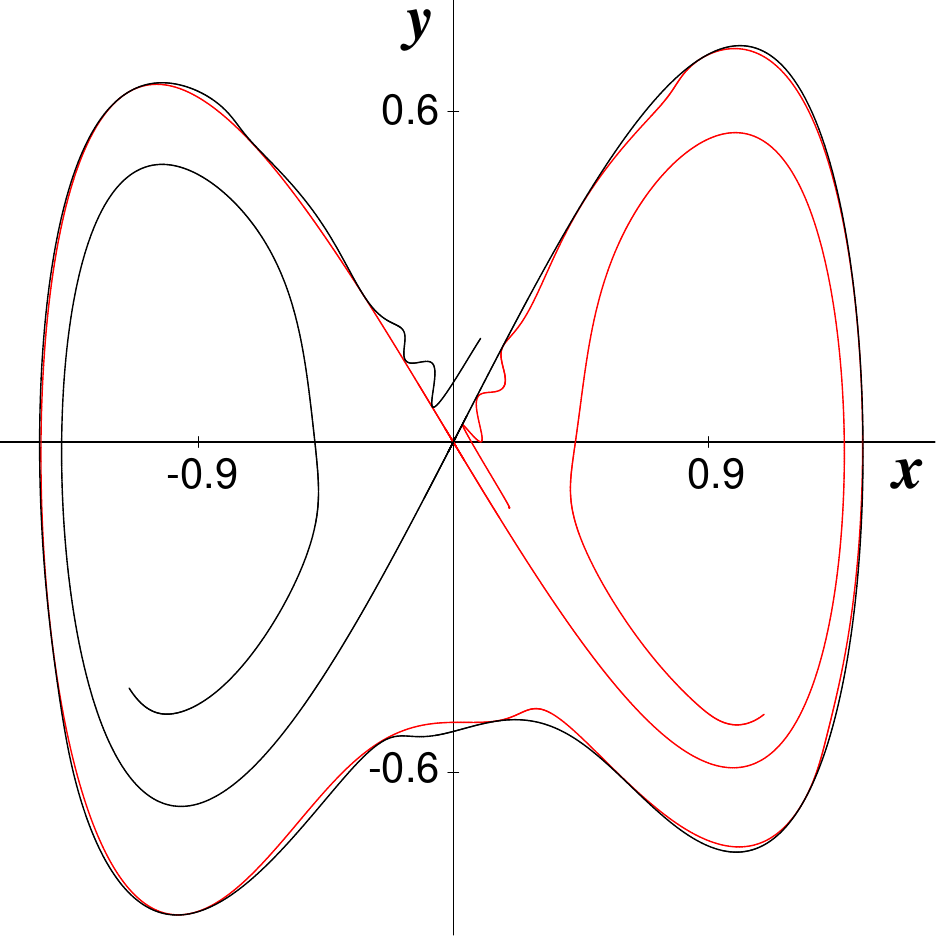} &
      \includegraphics[scale=0.4]{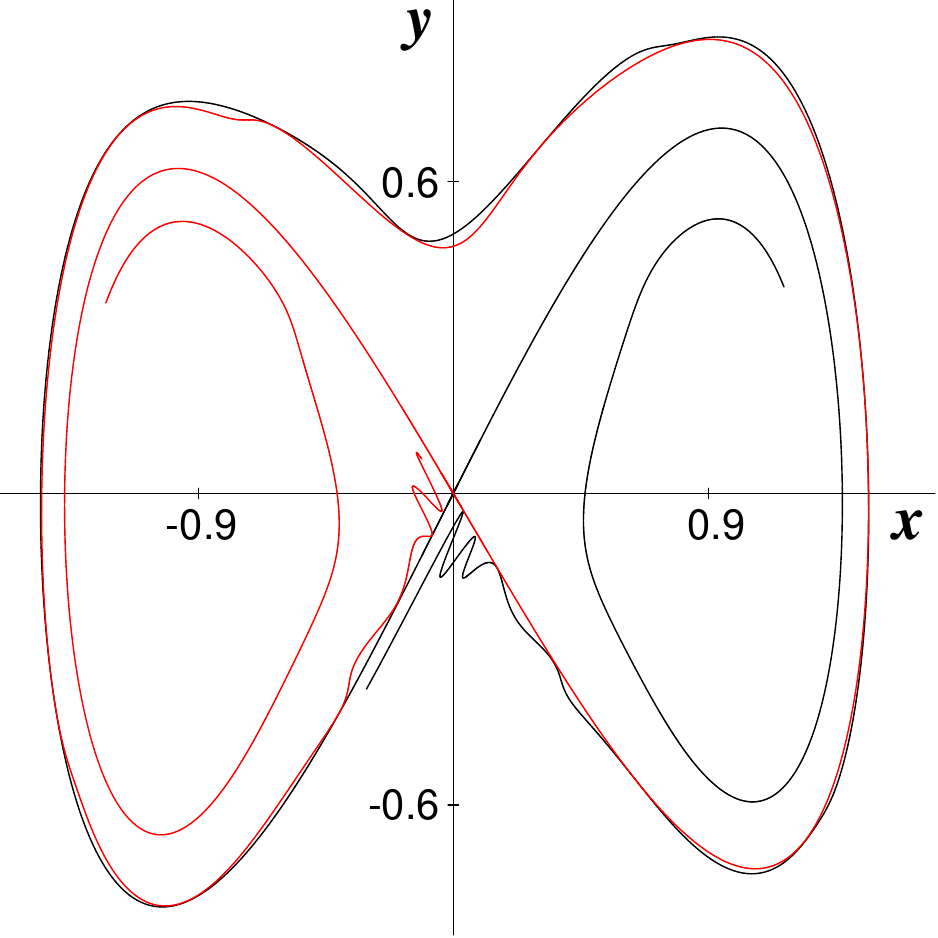}\\
      \footnotesize{(c)}& \footnotesize{(d)}\\
  \end{tabular}
\end{center}
\caption{ Behavior of separatrices of the fixed point $(0,0)$ for
Eq.~(\ref{eq31}) at  $\varepsilon=0.175$, $p_1=0.78549$,
$p_3=1.02$, $p_4=4$ and (a) $p_2=1.6$, (b) $p_2=-1.6$;
$\varepsilon=0.175$, $p_1=0.7850145$, $p_3=0.57$, $p_4=4$, (c)
$p_2=0.5$, (d) $p_2=-0.5$. } \label{fig10}
\end{figure}

The separatrices of the fixed point $(0,0)$ of the Poincar\'{e}
map on the $(\xi = x, \dot {\xi} = y) $ plane are shown in Fig.~
\ref{fig8} for  $\varepsilon=0.3$, $p_1=0.7551195621$,
$p_2=0.053875454$, $p_4=4$ and   $p_3=1.13$ (a), $p_3=1.7$ (b),
and  $p_3=2.83$ (c).

Note that for $p_1=0.8$, $p_2=0$, $p_3=0$  the unstable limit
cycle in Eq.~(\ref{eq1}) coincides with ``figure-eight''. Then,
for small enough $p_3\neq 0$, the inverse Poincar\'{e} map has a
quasiattractor.

When a perturbed autonomous equation has a ``big'' separatrix loop,
the Melnikov formula does not hold for a nonautonomous equation.
For this case, the separatrices of a fixed saddle point of the
Poincar\'{e} map for Eq.~(\ref{eq31}) on the $(\xi = x, \dot {\xi}
= y) $ plane are shown in Fig.~ \ref{fig9} for
 $ \varepsilon=0.1$, $p_1=0.78549$, $p_3 =1.02$, $p_4=4 $ and (a)
$p_2=1.6 $, (b) $p_2 =-1.6 $.  There occurs transversal
intersection of  the corresponding separatrices  ($ \Gamma
^{r}_{u}\bigcap \Gamma ^{l}_{s} \neq \oslash $ in
Fig.~\ref{fig9}(a) and $\Gamma ^{r} _ {s} \bigcap \Gamma ^{l}_{u}
\neq \oslash $ in Fig.~\ref{fig9}(b)). Also, homoclinic structures
with tangency (Fig.~\ref{fig10}) are possible at different values
of the parameters.

\begin{figure}[htb]
\begin{center}
\begin{tabular}{ccc}
\includegraphics[scale=0.4]{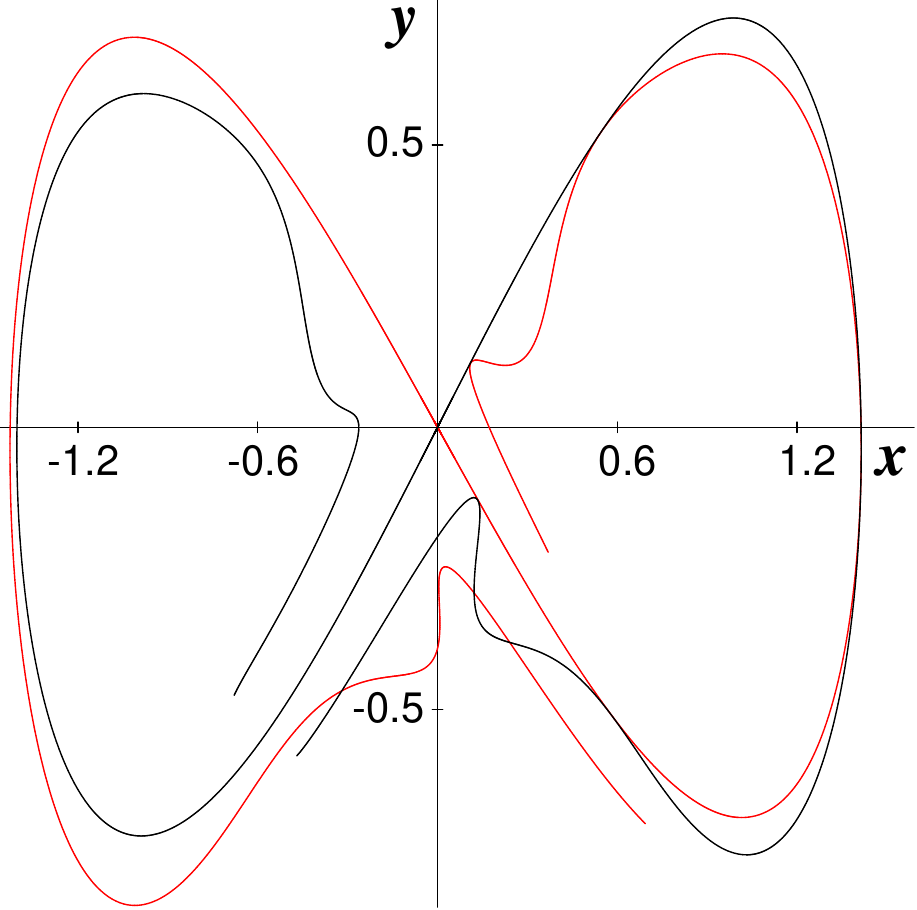} &
\includegraphics[scale=0.4]{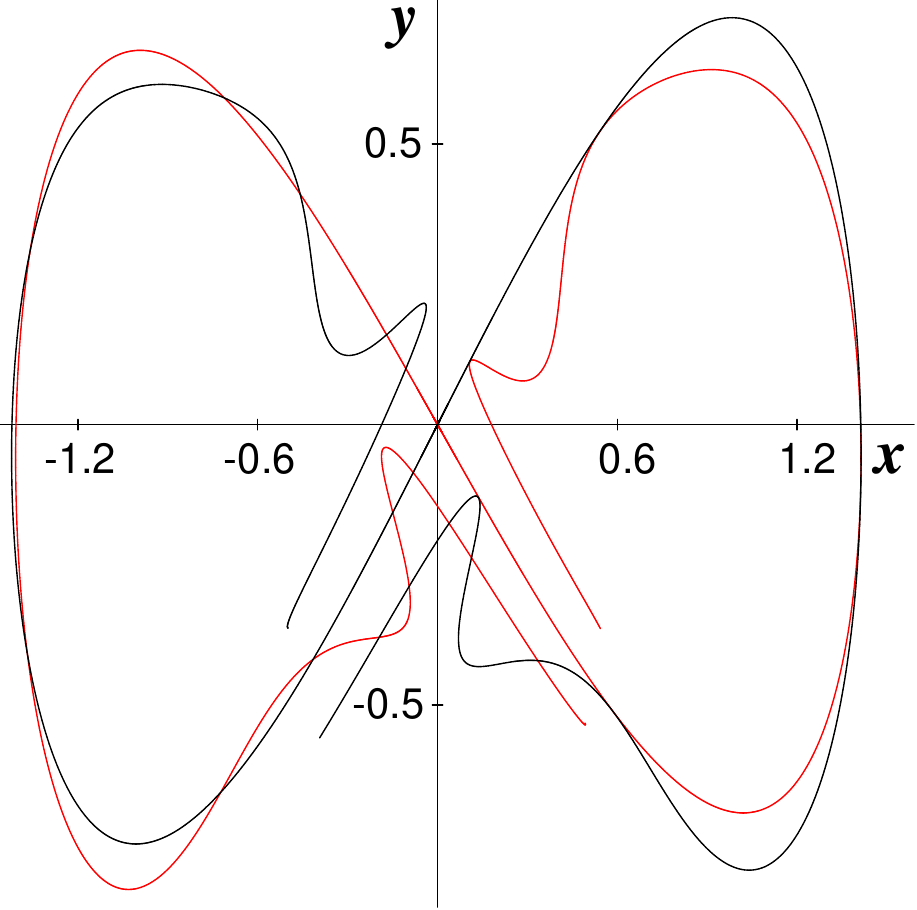}&
\includegraphics[scale=0.4]{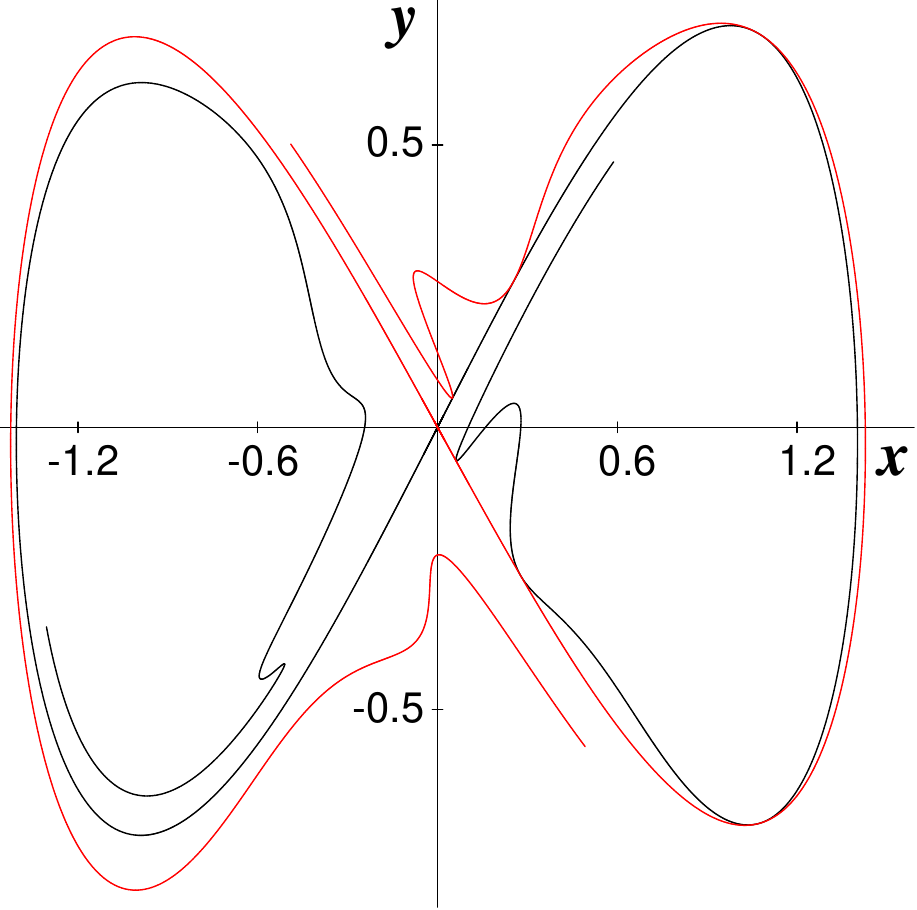}\\
\footnotesize{(a)}& \footnotesize{(b)} & \footnotesize{(c)}\\
\includegraphics[scale=0.4]{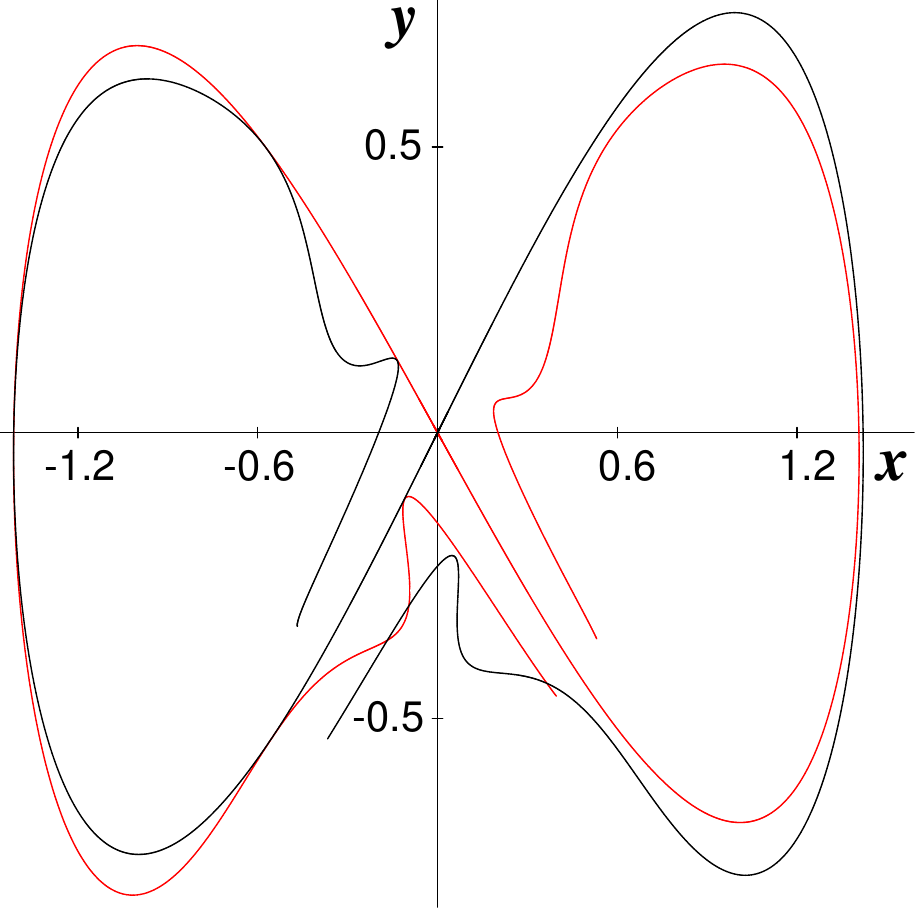}&
\includegraphics[scale=0.4]{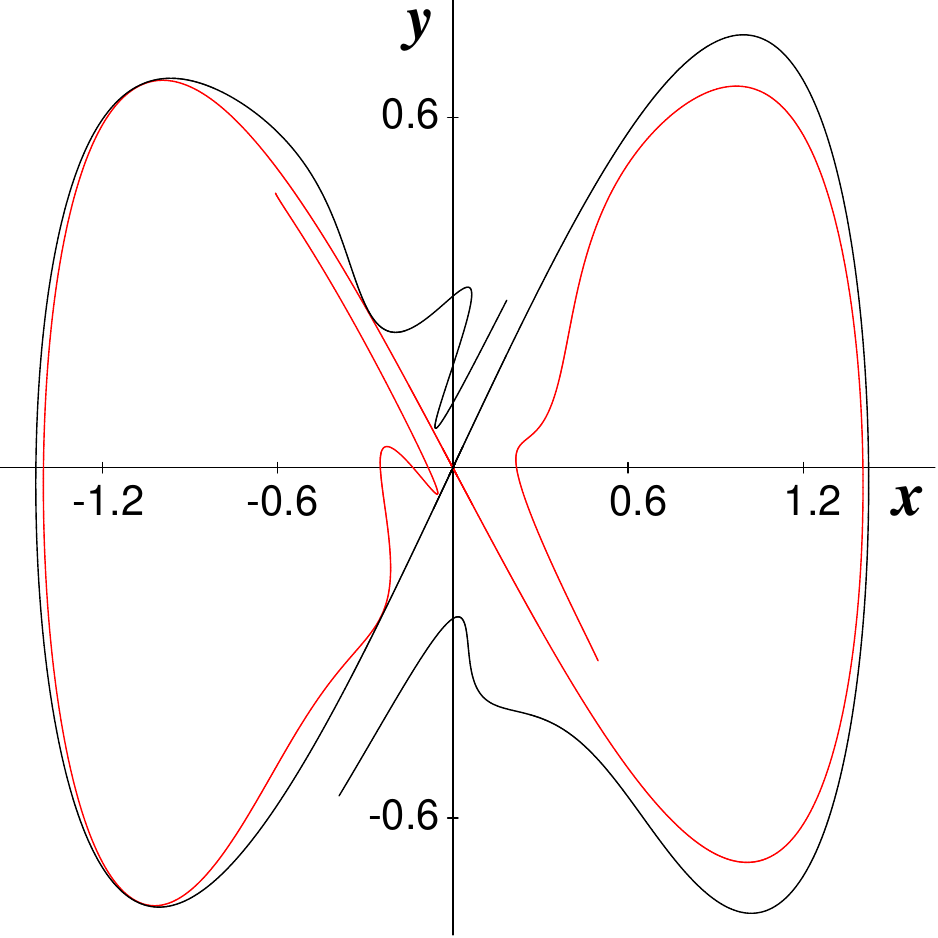} &
\includegraphics[scale=0.4]{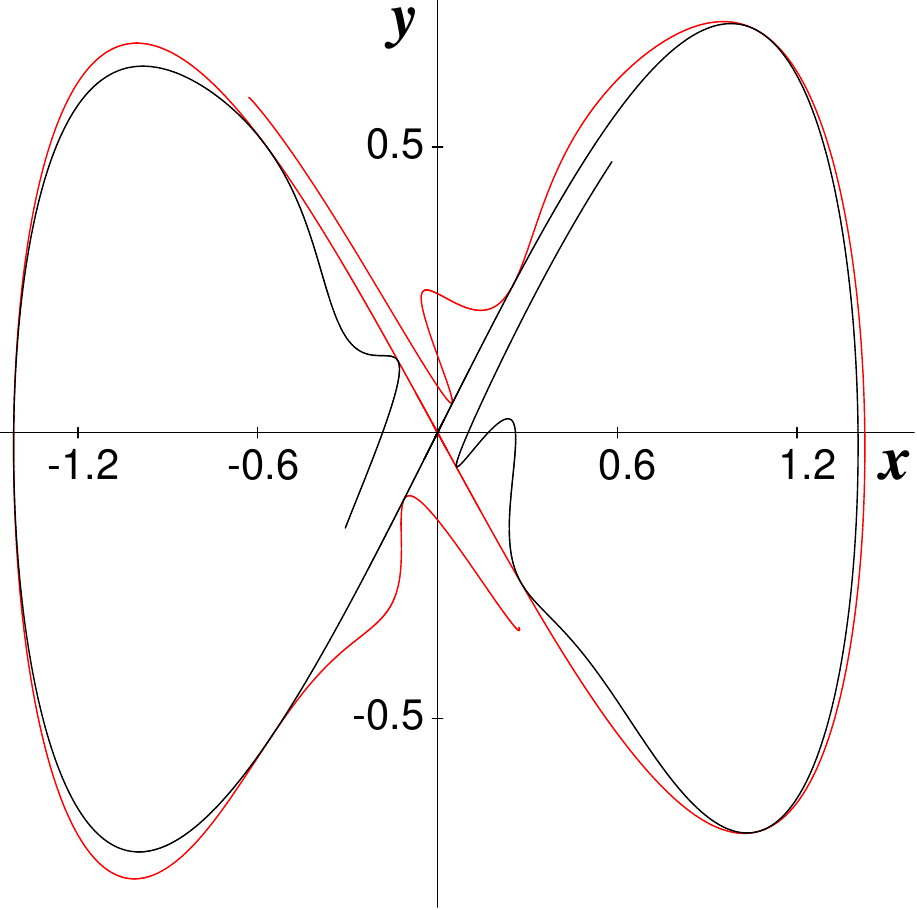}\\
 \footnotesize{(d)}&\footnotesize{(e)}& \footnotesize{(f)}\\
\includegraphics[scale=0.4]{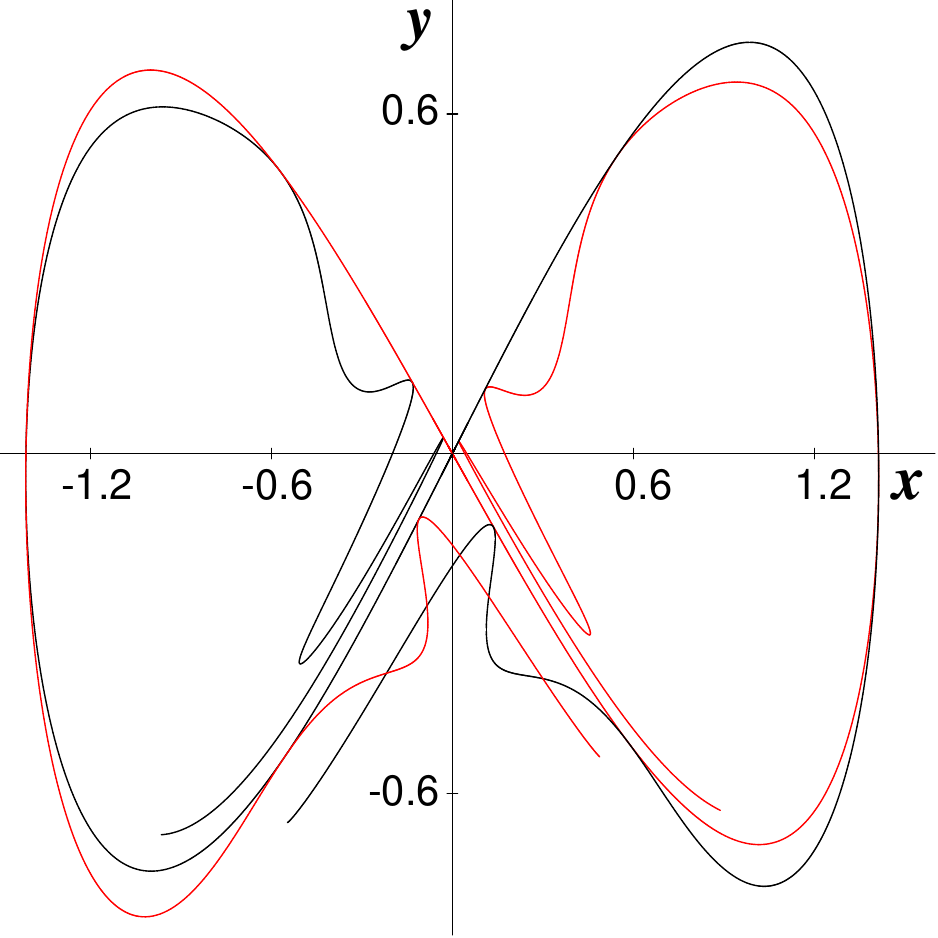} &
\includegraphics[scale=0.4]{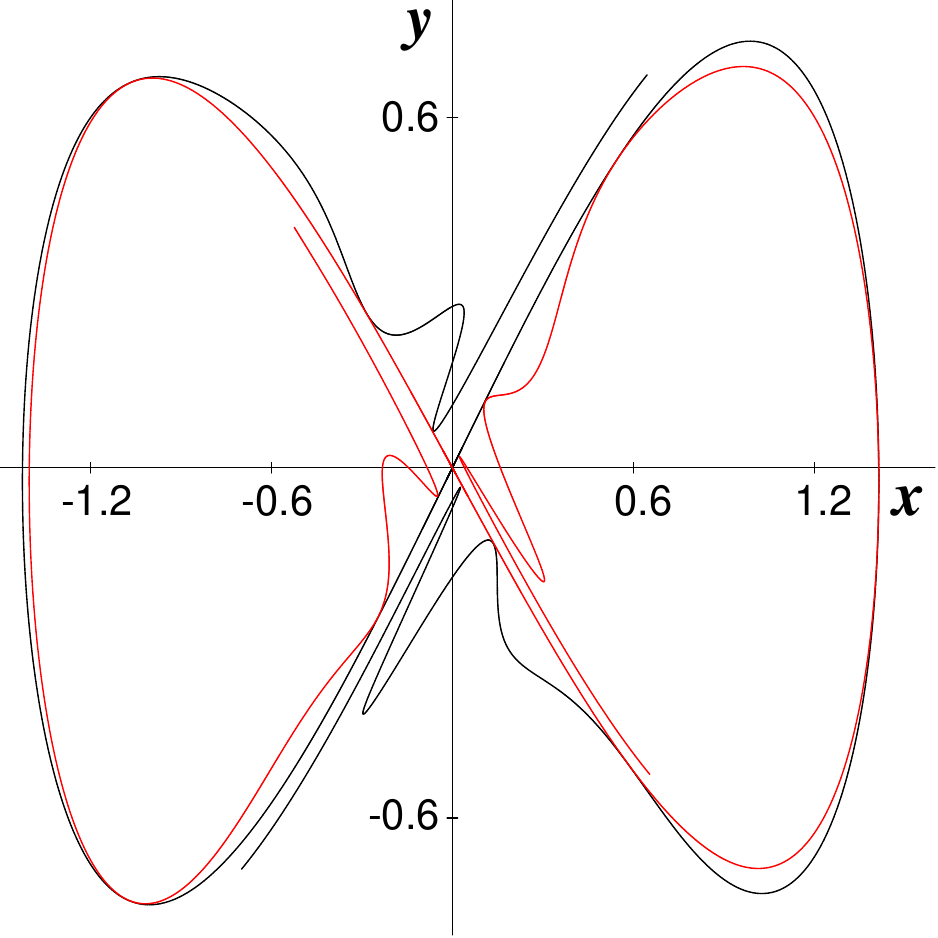}&
\includegraphics[scale=0.4]{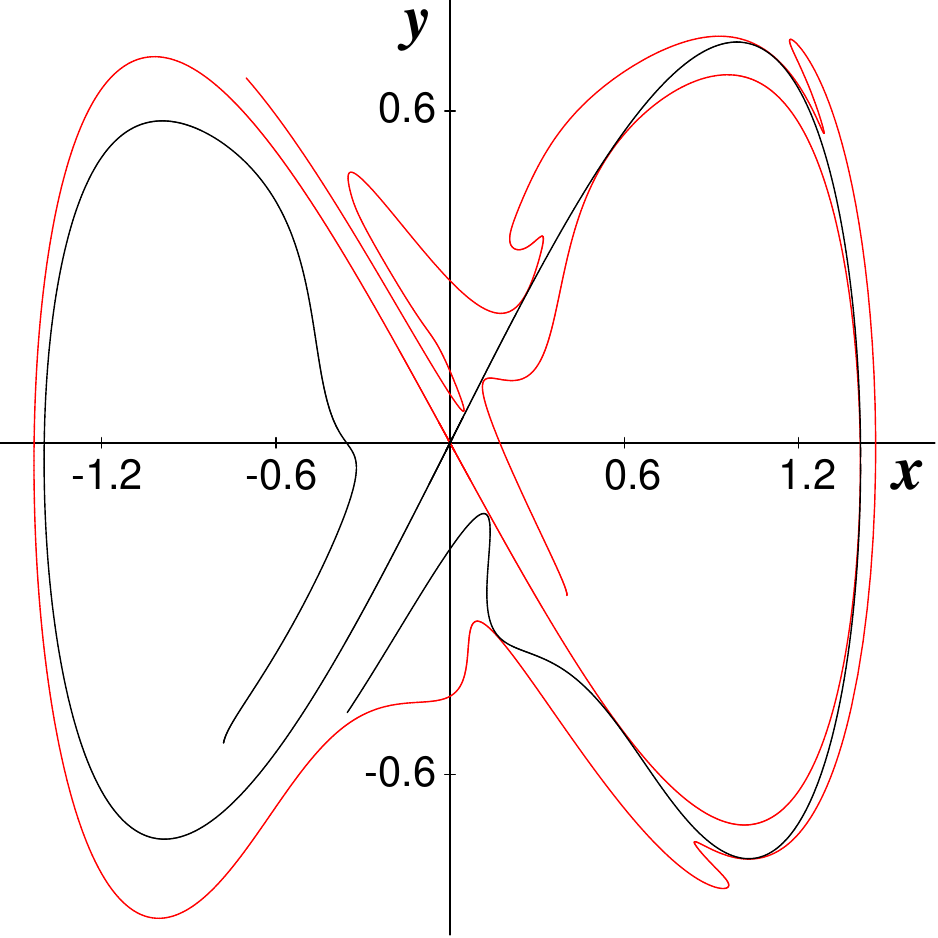} \\
\footnotesize{(g)}& \footnotesize{(h)}& \footnotesize{(i)}\\
\includegraphics[scale=0.4]{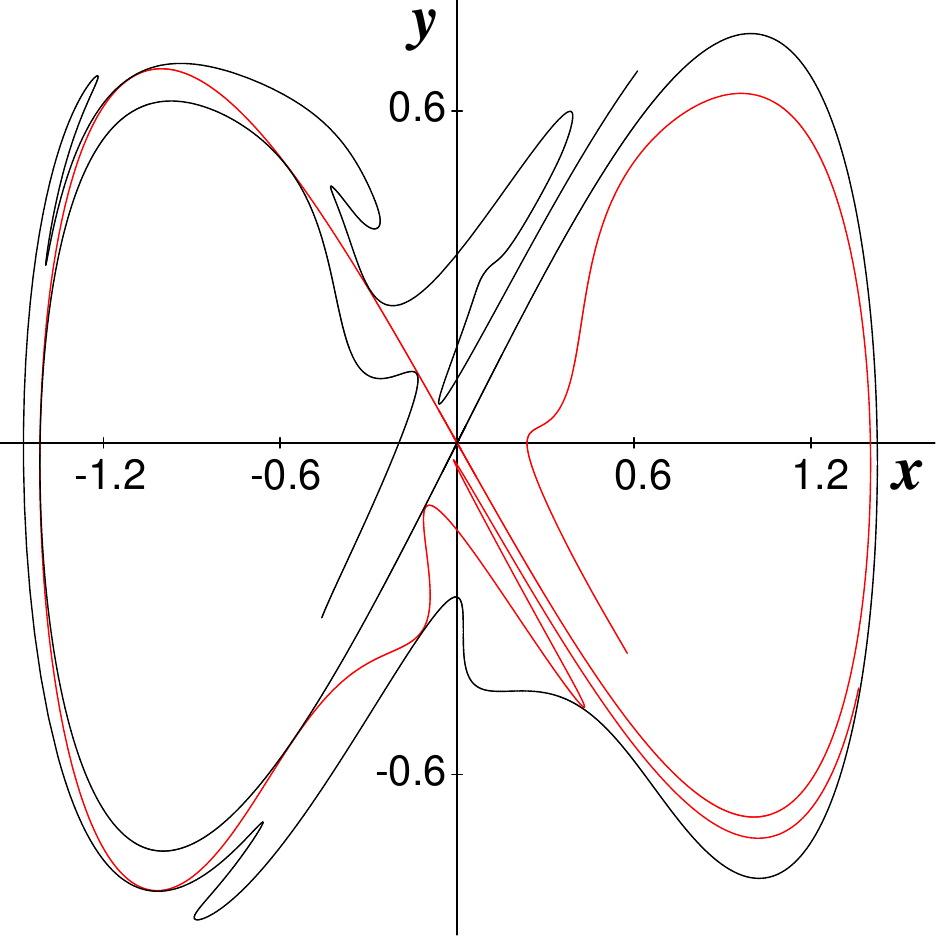}\\
 \footnotesize{(j)}\\
\end{tabular}
\end{center}
\caption{Other homoclinic structures with tangency of stable and
unstable separatrices of the fixed point $(0,0)$ for
Eq.~(\ref{eq31}).} \label{fig11}
\end{figure}

Figure \ref{fig11} illustrates other homoclinic structures  with
tangency of stable and unstable separatrices of the fixed point
$(0,0)$ for Eq.~(\ref{eq31}) at  $\varepsilon=0.12$, $p_4=4$ for
the following values of the parameters $p_1$, $p_2$, $p_3$: (a)
$p_1=0.7$, $p_2=0.3$, $p_3=3$; (b) $p_1=0.86$, $p_2=0.2$,
$p_3=4.55$; (c) $p_1=0.6$, $p_2=0.1$, $p_3=2.34$; (d) $p_1=0.86$,
$p_2=0.25$, $p_3=2.96$; (e) $p_1=1$, $p_2=0.1$, $p_3=2.32$; (f)
$p_1=0.7$, $p_2=0$, $p_3=2$; (g) $p_1=0.8$, $p_2=0.2$, $p_3=3.34$;
(h) $p_1=0.9$, $p_2=0$, $p_3=1.98$;  (i) $p_1=0.65$, $p_2=0.35$,
$p_3=2.82$; (j) $p_1=0.9$, $p_2=0.3$, $p_3=2.97$.

\section{Bifurcation diagrams}
Using the WInSet and Maple 13 software, we constructed three
bifurcation diagrams of the Poincar\'{e} map for Eq.~(\ref{eq31})
on the $(p_2, p_3)$ plane for fixed values of the parameters
$\varepsilon$, $p_1$, and $p_4$. In the bifurcation curves, the
corresponding separatrices of the fixed point  (0,0) are tangent
to each other. These curves separate domains with homoclinic
structure on the plane of parameters $(p_2, p_3)$ (the stable and
unstable separatrices of the saddle point (0,0) intersect
transversally). The three bifurcation diagrams describe all
possible cases of the  relative position of separatrices of the
fixed saddle point (0,0) for the Poincar\'{e} map.

\begin{figure}[htb]
\begin{center}
\includegraphics[scale=0.85]{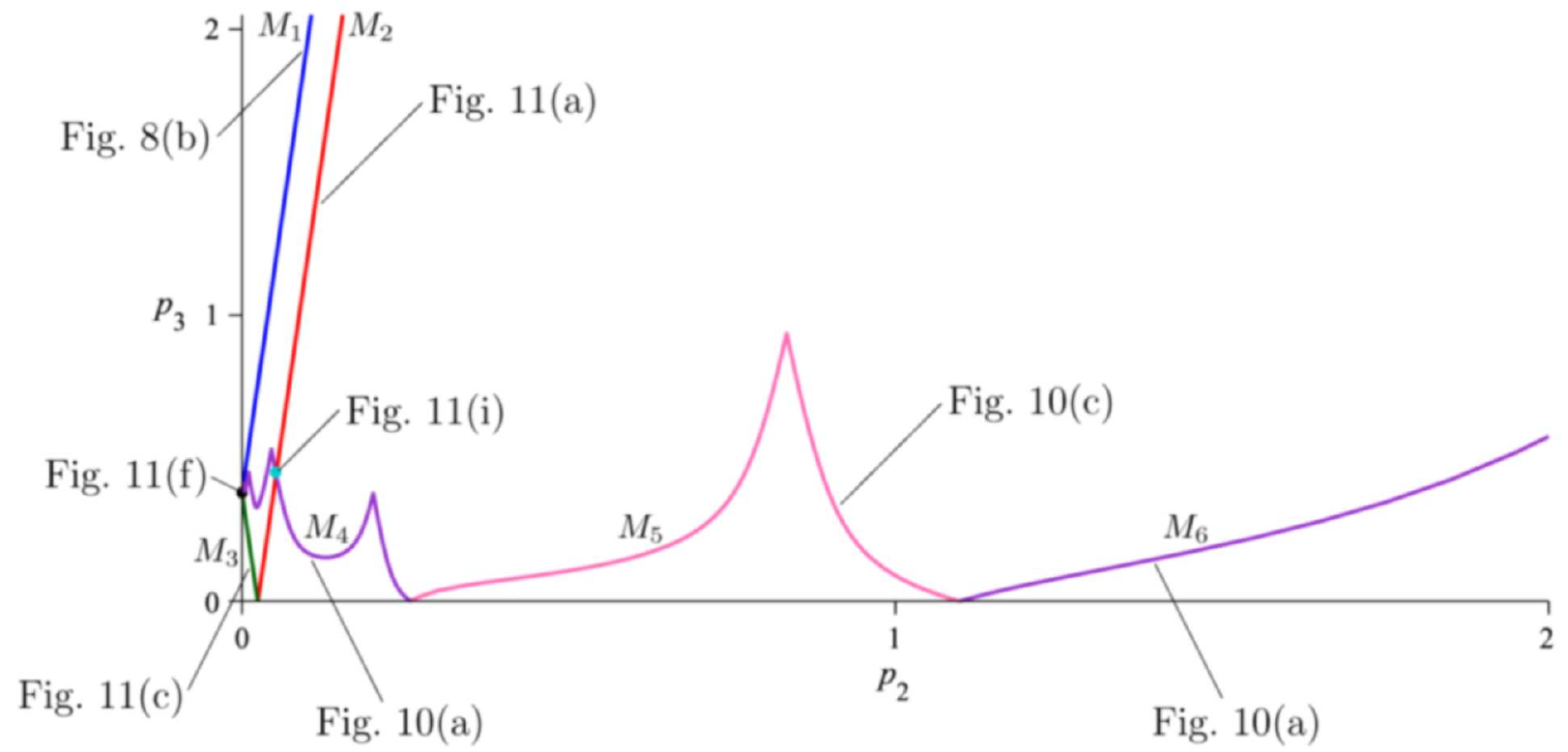}
\end{center}
\caption{Bifurcation diagram for the Poincar\'{e} map on the
$(p_2,p_3)$ plane at $p_1=0.78$.}\label{fig12}
\end{figure}

The obtained bifurcation diagrams are symmetric to the $p_3$ axis.
Let us set  $p_2>0$ and consider in more detail each of the three
bifurcation diagrams. By fixing $\varepsilon=0.12$, $p_1=0.78$,
$p_4=4$, we obtain six bifurcation curves. Equations for the
straight lines $M_1$, $M_2$, $M_3$ are found from (\ref{eq34}).
The other bifurcation curves $M_4$, $M_5$, $M_6$ are obtained
numerically by means of the WInSet software. Each pair of lines
$M_2$ and $M_3$, $M_4$ and $M_5$, $M_5$ and $M_6$ have exactly one
common point on the $p_2$ axis. The first point ($p_2\approx
0.024$) corresponds to the right separatrix loop in the autonomous
equation; the next two points ($p_2\approx 0.25838$ and
$p_2\approx 1.0983$) correspond to the ``big'' separatrix loop in
the autonomous equation. The intersection points of the lines
$M_2$ and $M_4$, as well as of the lines $M_1$, $M_3$ and $M_4$
correspond to double homoclinic tangency. The obtained bifurcation
curves are presented in Fig.~\ref{fig12}.

\begin{figure}[htb]
\begin{center}
\includegraphics[scale=0.85]{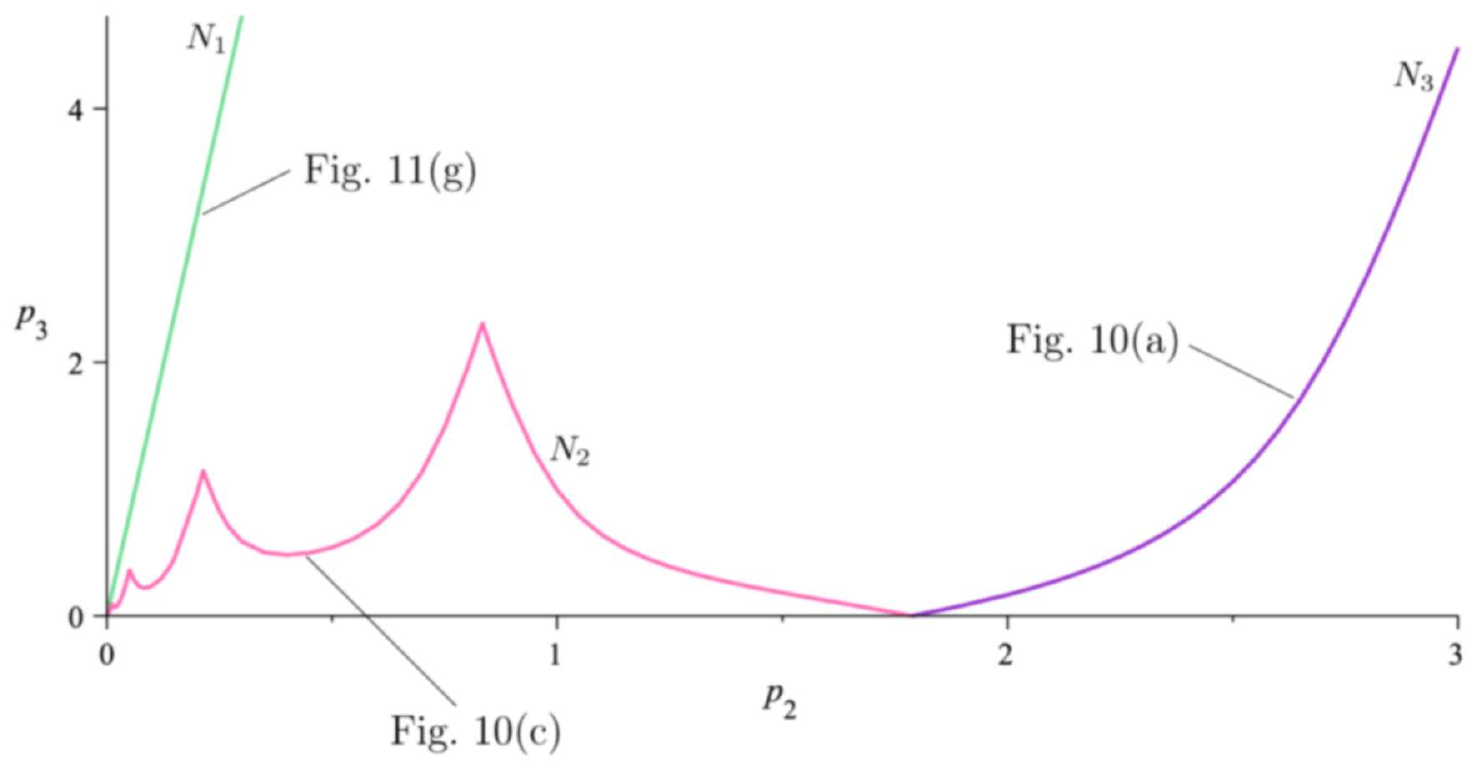}
\end{center}
\caption{Bifurcation diagram for the Poincar\'{e} map on the
$(p_2,p_3)$ plane at $p_1=0.8$.}\label{fig13}
\end{figure}

Setting $\varepsilon=0.12$, $p_1=0.8$, $p_4=4$, we obtain three
bifurcation curves shown in Fig.~\ref{fig13}. As the parameter
$p_1$ changes from $0.78$ to $0.8$ the straight lines $M_1$ and
$M_2$ in Fig.~ \ref{fig12} approach each other and coincide at
$p_1=0.8$. As a result we obtain a straight line $N_1$ with a new
type of tangency -- double homoclinic tangency. An equation for
$N_1$ is found from (\ref{eq34}). The lines $N_2$ and $N_3$
obtained numerically have one common point ($p_2\approx 1.788$,
$p_3=0$) corresponding to the ``big'' separatrix loop in the
autonomous equation.

\begin{figure}[htb]
\begin{center}
\includegraphics[scale=0.85]{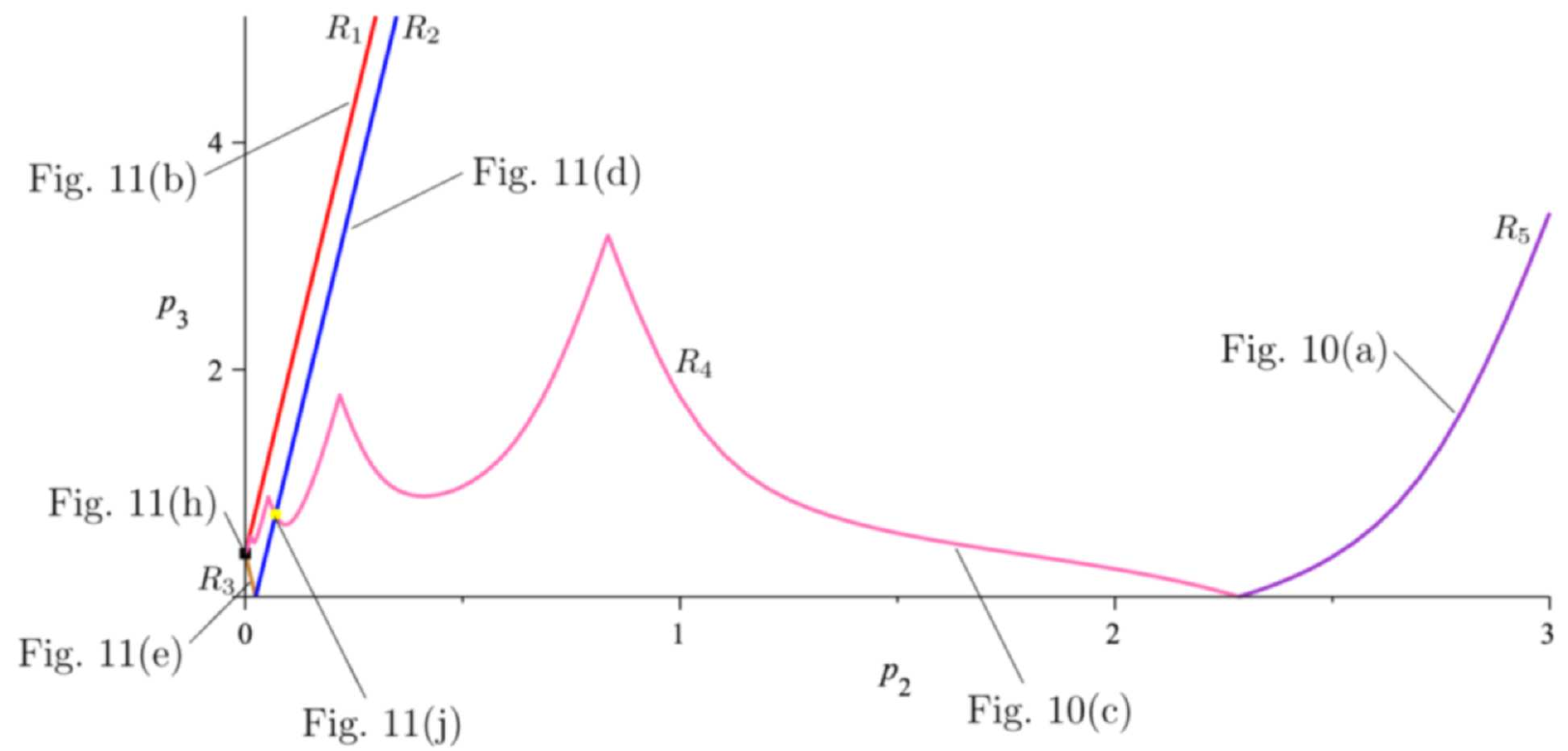}
\end{center}
\caption{ Bifurcation diagram for the Poincar\'{e} map on the
$(p_2,p_3)$ plane at $p_1=0.82$.}\label{fig14}
\end{figure}

Setting $\varepsilon=0.12$, $p_1=0.82$, $p_4=4$, we obtain five
bifurcation curves plotted in Fig.~\ref{fig14}. Equations for the
straight lines $R_1, R_2, R_3$ are found from (\ref{eq34}). The
other bifurcation lines  $R_4, R_5$ are obtained numerically using
the WInSet software. The intersection point of the curves $R_4$
and $R_5$ ($p_2\approx 2.28515$, $p_3=0$) corresponds to the ``big''
separatrix loop in the autonomous equation. The intersection
points of $R_2$ and $R_4$ and of $R_1$, $R_3$ and $R_4$ give
double homoclinic tangencies.

Each of the three bifurcation diagrams has domains with homoclinic
structure and a nonsmooth boundary. This phenomenon was explained
in ample detail in the work \cite {GSV}.

\section{Conclusion}
The problem of time-periodic perturbations of two-dimensional
Hamiltonian systems with a saddle and two separatrix loops in the
form of ``figure-eight'' is a challenging problem for the theory of
bifurcations. Bifurcations in the neighborhood of ``figure-eight'' for
the case of an unperturbed autonomous system with a nonzero saddle
value were recently considered in \cite {GSV}. This problem for
the case of a zero saddle value has not been fully understood yet.
The asymmetric Duffing--Van-der-Pol equation (\ref {eq1}) studied
in the present work is a good model for solution of this problem.

Despite its fundamental role in the theory of differential
equations, the theory of bifurcations, and the theory of
oscillations, Eq.~(\ref {eq1}) has not been studied thus far for
the case when $p_2\ne 0$.  We have solved the problem of limit
cycles in the autonomous case. For the nonautonomous case, we have
found resonance zone structures and global behavior of solutions
in the cells separated from unperturbed separatrices. Different
resonance periodic solutions and two-dimensional invariant tori
have also been found. The problem of the existence of homoclinic
structures in the neighborhood of unperturbed separatrices (in the
neighborhood of ``figure-eight'') has been solved. All possible cases
of relative position of the separatrices of a trivial fixed saddle
point  for the Poincar\'{e} map have been revealed. Three
bifurcation diagrams for the Poincar\'{e} map on the $(p_2,p_3)$
plane separating domains of existence of different homoclinic
structures have been constructed. The results obtained for the
separatrix tangency illustrate many specific features found in
\cite {GSV} for two-parametric families of maps in the
neighborhood of ``figure-eight'' with nonzero saddle value.

Note that some of the problems associated with the presence of homoclinic structures remain open. For example, one such problem is to study fractal properties of attraction basin boundaries for stable periodic regimes in the considered equation. 

\section{Acknowledgments}

We dedicate this paper to the memory of the outstanding scientist Leonid P. Shil'nikov, the pioneer of homoclinic bifurcation theory.
The authors are grateful to M.I. Malkin for helpful discussions
and comments.
This work was partially supported by the Russian Science Foundation, grant No 14-41-0044.

\end{document}